\documentstyle[mssymb]{article}
\def\be{\begin{equation}}
\def\ee{\end{equation}}
\def\bea{\begin{eqnarray}}
\def\eea{\end{eqnarray}}
\def\Qed{\rule{1ex}{1ex}}
\def\B{{\sf B}}
\def\R{{\bf R}}
\def\Z{{\bf Z}}
\def\DD{{\cal D}}
\def\FF{{\cal F}}
\def\SS{{\cal S}}
\def\eps{\epsilon}
\def\leqs{\leqslant}
\def\geqs{\geqslant}
\def\ra{\rightarrow}
\def\iff{\Leftrightarrow}
\def\nl{\hspace*{\fill}\\}
\def\nls{\nl\hspace*{2em}}
\def\ul{\underline}
\def\romit#1{\item[{\rm #1}]}
\def\x{\makebox[.4cm]{$x$}}
\def\b{\makebox[.4cm]{$\bullet$}}
\def\0{\makebox[.4cm]{$\circ$}}
\def\:{\makebox[.4cm]{$\cdot$}}
\def\1#1{\makebox[.4cm]{$#1$}}
\def\sq#1{\begin{array}{c}#1\end{array}}
\begin{document}
\centerline{\large
The still-Life density problem and its generalizations}
\vspace*{2ex}
\centerline{Noam D.\ Elkies}
\centerline{June 1997}

\vspace*{5ex}

\renewcommand{\baselinestretch}{.9}
{\small
{\bf Abstract.} A {\em still Life} is a subset~$S$\/ of the square
lattice $\Z^2$ fixed under the transition rule of Conway's Game of
Life, i.e.\ a subset satisfying the following three conditions:
\begin{description}
\romit{\qquad 1.} No element of $\Z^2-S$ has exactly three neighbors
 in~$S$\/;
\romit{\qquad 2.} Every element of~$S$\/ has at least two neighbors
 in~$S$\/;
\romit{\qquad 3.} Every element of~$S$\/ has at most three neighbors
 in~$S$.
\end{description}
Here a ``neighbor'' of any $x\in\Z^2$ is one of the eight lattice points
closest to~$x$ other than~$x$ itself.  The {\em still-Life conjecture}
is the assertion that a still Life cannot have density greater than~1/2
(a bound easily attained, for instance by $\{ x: 2|x_1 \}$).  We prove
this conjecture, showing that in fact condition~3 alone ensures that
$S$\/ has density at most~1/2.  We then consider variations of the
problem such as changing the number of allowed neighbors or the
definition of neighborhoods; using a variety of methods we find
some partial results and many new open problems and conjectures.
}
\renewcommand{\baselinestretch}{1}

\vspace*{5ex}

{\bf 1.~Life, the still-Life conjecture, Conjecture~A, and $\delta(n)$.}

Let $L$ be the square lattice $\Z^2$ in the plane.  Two points
$x,y$ of~$L$ are said to be {\em neighbors}, or {\em adjacent},
if $x\neq y$ but $|x_i-y_i|\leqs1$ for $i=1,2$; equivalently,
if $x,y$ are at distance $1$ or~$\sqrt2$.  Thus each $x\in L$ has
8 neighbors.  For any subset $S\subseteq L$ and $x\in L$ let
$N_S(x)$ be the number of neighbors of~$x$ contained in~$S$,
which is an integer in $[0,8]$.  The {\em lower} and {\em upper
densities} of $S\subseteq L$ are defined as the $\liminf_{r\ra\infty}$
and $\limsup_{r\ra\infty}$ of $|\B_r\cap S| / |\B_r|$, where $\B_r$
is the box $\{ x \in L: |x_1|, |x_2| < r \}$ and $|\cdot|$ is
used for the cardinality of a set.  These are real
numbers in $[0,1]$; if they are equal, their common value is
called the {\em density} of~$S$.  (This allows the set
$\{ x\in L: x_1>0 \}$ to have density $1/2$, a kind of pathology
one can circumvent by formulating a more restrictive definition;
but our definition is sufficient for the purposes at hand.)
If there is a subgroup $L'\subseteq L$ of finite index such
that $S$\/ is invariant under translation by~$L'$ we say
$S$\/ is {\em periodic}, and call the largest such $L'$ the
{\em period lattice} of~$S$.  A periodic set necessarily has
a density which is a rational number (with denominator a factor
of the index in~$L$ of the period lattice).
We will denote the set $2^L$ of all subsets $S\subseteq L$ by~$\SS$.

Conway's Game of Life \cite[Ch.25]{WW} is in effect the iteration 
of a map $\lambda: \SS \ra \SS$ defined as follows: $x\in\lambda S$\/
$\iff$ either $x\in S$ and $2\leqs N_S(x) \leqs 3$, or $x\notin S$ and
$N_S(x)=3$.  For instance, $\lambda(\emptyset)=\lambda(L)=\emptyset$.
A subset $S\in L$ fixed by $\lambda$ is called a {\em still Life}.
Thus $S$\/ is a still Life if and only if:
\begin{description}
\romit{\qquad 1.} $x\notin S \Rightarrow N_S(x) \neq 3$;
\romit{\qquad 2.} $x\in S \Rightarrow N_S(x) \geqs 2$;
\romit{\qquad 3.} $x\in S \Rightarrow N_S(x) \leqs 3$.
\end{description}
For instance, we just saw that $\emptyset$ is a still Life.
So are the ``block'' $\{0,1\} \times \{0,1\}$,
which is the smallest nonempty still Life, and many infinite subsets
such as $\{ x: x_1 \equiv 0 {\rm\ or\ 1} \bmod m,\;
x_2 \equiv 0 {\rm\ or\ 1} \bmod n \}$ for any $m,n\geqs3$
(a lattice of blocks), or $\{ x: 2|x_1 \}$;  see Figure~1.
[We use $\bullet$ for an element of $S$, and $\circ$ or $\cdot$
 for an element of $L-S$; an empty space or one marked ``?''
 means a point which may be in either $S$ or $L-S$.]
It has long been noticed that while there are still-Life patterns
of density $\leqs 1/2$, including many that attain density~1/2
in different ways (see Appendix~A), no still Life was found to exceed
that density.  This gave rise to the {\bf still-Life conjecture}:
the density of any still Life is at most~$1/2$.

$$
\sq{
\:\:\:\:\:\:\:\:\:\:\:\:\\
\:\b\b\:\:\b\b\:\:\b\b\:\\
\:\b\b\:\:\b\b\:\:\b\b\:\\
\:\:\:\:\:\:\:\:\:\:\:\:\\
\:\:\:\:\:\:\:\:\:\:\:\:\\
\:\:\:\:\:\:\:\:\:\:\:\:\\
\:\b\b\:\:\b\b\:\:\b\b\:\\
\:\b\b\:\:\b\b\:\:\b\b\:\\
\:\:\:\:\:\:\:\:\:\:\:\:\\
\:\:\:\:\:\:\:\:\:\:\:\:\\
\:\:\:\:\:\:\:\:\:\:\:\:\\
\:\b\b\:\:\b\b\:\:\b\b\:\\
\:\b\b\:\:\b\b\:\:\b\b\:\\
\:\:\:\:\:\:\:\:\:\:\:\:\\
\:\:\:\:\:\:\:\:\:\:\:\:
}
\qquad
\sq{
\:\b\:\b\:\b\:\b\:\b\:\\
\:\b\:\b\:\b\:\b\:\b\:\\
\:\b\:\b\:\b\:\b\:\b\:\\
\:\b\:\b\:\b\:\b\:\b\:\\
\:\b\:\b\:\b\:\b\:\b\:\\
\:\b\:\b\:\b\:\b\:\b\:\\
\:\b\:\b\:\b\:\b\:\b\:\\
\:\b\:\b\:\b\:\b\:\b\:\\
\:\b\:\b\:\b\:\b\:\b\:\\
\:\b\:\b\:\b\:\b\:\b\:\\
\:\b\:\b\:\b\:\b\:\b\:\\
\:\b\:\b\:\b\:\b\:\b\:\\
\:\b\:\b\:\b\:\b\:\b\:\\
\:\b\:\b\:\b\:\b\:\b\:\\
\:\b\:\b\:\b\:\b\:\b\:\\
}
$$
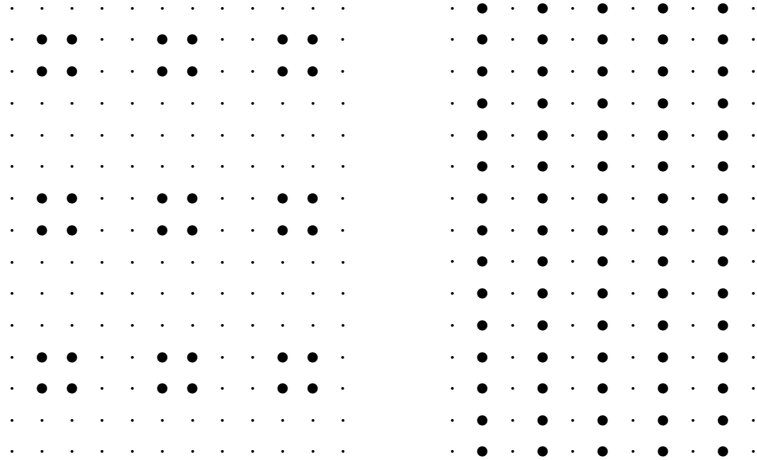
\centerline{Figure 1: Some still Lifes}

Trying to refute this conjecture by constructing a denser pattern~$S$,
one quickly finds that only condition (3) on~$S$\/ causes any trouble:
in any periodic pattern that comes close to attaining, let alone
exceeding, density 1/2, no $x\notin S$ has $N_S(x)<4$ (and even
$N_S(x)=4$ is rare), and no $x\in S$ has $N_S(x)<2$.  That is,
the following stronger conjecture is suggested:

{\bf Conjecture A.}  {\sl If, for some $S\in L$, every $x\in S$
has at most 3 neighbors in~$L$, then $S$ has (upper) density 
at most $1/2$ in~$L$.}

More tersely, if we define for any $S\in L$ its maximum degree
$d(S)$ by
$$
d(S) = \max_{x\in S} N_S(x)
$$
then Conjecture~A asserts that if $d(S)\leqs3$ then
$S$\/ has (upper) density at most $1/2$ in~$L$.
We put ``upper'' in parentheses because (the condition
$d(S)\leqs 3$ being local) from a subset $S\in L$
with upper density $\delta$ we readily construct
$S'$ of density~$\delta$ with the same maximum degree.
To see this, first note that for fixed $r_0\geqs1$ we can
omit $O(r)$ points from $\B_r$ and tile the rest with
translates of $\B_{r_0}$.   Thus if for each $\eps$
we may find arbitrarily large $r$ such that
$|\B_r\cap S|\geqs(\delta-\eps) |\B_r|$ then
at least one of the translates of $\B_{r_0}$ in our
tiling of $\B_r$ meets~$S$\/ in at least
$(\delta-\eps-O(1/r)) |\B_{r_0}|$ points
(since this is true on average over all those translates).
Since $\eps$ and $1/r$ are arbitrarily small it follows
that there is a subset $S_{r_0} \subset \B_{r_0}$ of size at
least $\delta \cdot |\B_{r_0}|$ with $d(S_{r_0})\leqs3$.
Removing some points from $S_{r_0}$ (which will not increase
its maximum degree), we may assume that
$|S_{r_0}| = \lfloor \delta \cdot |\B_{r_0}| \rfloor $.
Now, for $m=2,3,\ldots$, tile the region $\{x\in\R^2 :
4^{m-1} < \max(|x_1|,|x_2|) < 4^m\}$ with squares of
side-length $2^{m+1}$, and place a translate of $S_{2^m}\subset\B_{2^m}$
in each of these squares.  Note that no two of these
translates interact because they are distance at least~2 from
each other.  Thus their union is a set $S'\subset L$ with
$d(S)\leqs3$, and $S'$ has density~$\delta$ as claimed.

Note that Conjecture~A is mathematically much more natural
than the original still-Life conjecture (which Conj.~A implies),
and suggests various generalizations.  Thus for each integer~$n\geqs0$
we may ask what is the maximal density $\delta(n)$
of a subset $S\in L$ with $d(S)\leqs n$.  (We could generalize
further, to other lattices $L$ and beyond; we consider some of
these generalizations at the end of this paper.)  We might have
asked instead for the supremum of the upper densities of
such~$S$, but the same cut-and-paste argument shows that
this supremum is actually attained by a pattern with that
density.  Clearly $\delta(n)$ is a non-decreasing function
of~$n$, and we might expect it to increase strictly with
$n\in [0,8]$.  But in fact the example of $\{ x \in L: 2|x_1 \}$
shows that $\delta(2)\geqs 1/2$, so Conj.~A implies the surprising
$\delta(2)=\delta(3)$: when we allow maximum degree~3 rather than~2 
we permit a much greater variety of patterns attaining density~$1/2$
(as seen in Appendix~A, whereas we shall see that the $\delta(2)=1/2$
pattern is essentially unique), but none with density any higher
than~$1/2$.

To obtain a lower bound $\delta$ on $\delta(n)$ it is enough to exhibit
one $S\in\SS$ of density~$\delta$ with $d(S)\leqs n$.  We collect in
the next lemma the best such bounds known to us.

{\bf Lemma 1.}  {\sl For each $n\geqs0$, the corresponding $\delta$
listed in the following table is a lower bound on $\delta(n)$:
$$
\begin{array}{c|ccccccccr}
    n  &  0  &  1  &  2  &  3  &  4  &   5  &  6  &  7  & \geqs8
\\ \hline
\delta & 1/4 & 1/3 & 1/2 & 1/2 & 3/5 & 9/13 & 4/5 & 8/9 & 1
\end{array}
$$
}

{\em Proof}\/: For each $n$ we exhibit $S\in\SS$ with $d(S)\leqs n$
attaining the claimed $\delta$; see Appendix~B for the pictures.
For $n=0$ take
$$
S = \{ x \in L : x_1 \equiv x_2 \equiv 0 \bmod 2 \};
$$
for $n=1$,
$$
S = \{ x \in L : x_1 \equiv \pm1 \bmod 3, x_2 \equiv 0 \bmod 2 \}.
$$
For $n=2,\,3$ we noted already that $\{ x \in L: 2|x_1 \}$
has maximum degree~2 and density~$1/2$.  For $n=4$ we may take
$$
S = \{ x \in L : x_1+2x_2 \not\equiv 0,\,1 \bmod 5 \}
$$
(or the alternative period set of maximum degree~4, density $3/5$,
and period lattice $2\Z \times 5\Z$ also pictured in Appendix~B);
for $n=5$,
$$
S = \{ x \in L : 2x_1+3x_2 \not\equiv \pm 1, \; \pm 5 \bmod 13 \}
$$
(the four excluded values being the cubic residues mod~13).
The lattice complements
$$
L - \{ x : x_1+2x_2 \equiv 0 \bmod 5\}, \quad
L - \{ x : x_1 \equiv x_2 \equiv 0 \bmod 3 \}
$$
deal with $n=6$ and $n=7$.  Finally for $n\geqs8$ we of course
take $S=L$.~~\Qed

Note that all the sets $S$\/ appearing here are periodic, and
most of them have square period lattices.  A square period
lattice suggests that $S$\/ can be described more succinctly
by identifying $L$ with the ring $\Z[i]$ of Gaussian integers.
For instance, the set of maximum degree~5 and density $9/13$
may be described as $\{ x \in \Z[i]: 2x \not\equiv \pm1, \; \pm i
\bmod 3 - 2i \}$, and the set of degree~6 and density $4/5$
is the complement in $\Z[i]$ of the ideal $(2-i)$.

\vspace*{3ex}

{\bf 2.~Easy $\delta(n)$ values; the Voronoi-cell method.}

Most of the lower bounds in Lemma 1 can be easily proved sharp.
For the cases $n\geqs6$ we show this by adapting a double-counting
method familiar from the combinatorics of finite graphs:

{\bf Proposition 1.}
{\sl For each $n$ we have $\delta(n)\leqs 8/(16-n)$.
Moreover a periodic set $S$ with $d(S)\leqs n$ has density $8/(16-n)$
if and only if $N_S(y)=8$ for each $y\notin S$ (i.e.\ no two elements
of~$L-S$ are adjacent) and $N_S(x)=n$ for each $x\in S$.}

{\em Proof}\/: For some large $r$ count the adjacent
pairs $(x,y)$ in $\B_r$ with $x\in S$, $y\notin S$.  
Since $N_S(x)\leqs n$ for each $x\in S$\/ there are at least
$(8-n) |\B_{r-1} \cap S| = (8-n-O(1/r)) |\B_r \cap S|$ such pairs.
On the other hand, since $N_S(y)\leqs 8$ for each $y\notin S$\/
the number of pairs cannot exceed
$8 |\B_r - S| = 8 |\B_r| - 8 |\B_r\cap S|$.
Thus $(16-n-O(1/r)) |\B_r \cap S| \leqs 8 |\B_r|$.
Dividing by $|\B_r|$ and letting $r\ra\infty$ we conclude that
$S$\/ has upper density at most $8/(16-n)$.   Moreover if
some periodic $S$\/ with $d(S)\leqs n$ has density $8/(16-n)$
then equality must hold at each step, so $N_S(y)=8$ for each
$y\notin S$\/ and $N_S(x)=n$ for each $x\in S$\/ as claimed;
conversely if these hold then the same double-counting argument
shows that $|\B_r \cap S| = (8/(16-n)) |\B_r| + O(r)$ so
$S$\/ has density $8/(16-n)$.~~\Qed

Since we have already exhibited $S$\/ of maximal degree $n=6$,~7,~8
with density $8/(16-n)=4/5,\ 8/9,\ 1$, we conclude:

{\bf Corollary.}  {\sl For $n=6,7,8$ we have $\delta(n)= 8/(16-n)$.
Moreover, the periodic sets $S$ of (maximal) degree~$6$ and density
$4/5$ are $\{ x \in \Z[i]: x \not\equiv c \bmod 2-i \}$ and
$\{ x \in \Z[i]: x \not\equiv c \bmod 2+i \}$ for some $c$,
and the periodic sets $S$ of (maximal) degree~$7$ and
density $8/9$ are the complements of the density-$1/9$ sets
$\{ x \in L: x_1 \equiv c \bmod 3, x_2 \equiv b(x_1) \bmod 3 \}$
and $\{ x \in L: x_2 \equiv c \bmod 3, x_1 \equiv b(x_2) \bmod 3 \}$
for some integer $b \bmod 3$ and some periodic function
$b: 3\Z+c \ra \Z/3$.}

{\em Proof}\/: For $n=6,7,8$ our lower and upper bounds of Lemma~1
and Prop.~1 agree, so we need only verify that our list of periodic sets
attaining this bound for $n=6,7$ is complete.  For $n=7$ the condition
of Prop.~1 is equivalent to the requirement that $3\times 3$ squares
centered at the points of $L-S$ form a periodic tiling of $\R^2$,
which is readily seen to happen only if $L-S$ is one of the sets
of density $1/9$ exhibited above.  For $n=6$, let $S\subset L$\/
be any subset, periodic or not, each of whose elements has two
neighbors outside~$S$\/ while no two elements of $L-S$ are adjacent.
Fix some $x\in S$\/; without loss of generality (i.e.\ by applying
an isometry of~$L$) we may assume that $x=(0,0)$ and that its
two neighbors in $L-S$ are $(\pm1,\pm1)$ or $(\pm1,0)$ or $(\pm1,1)$
or $(1,0),(-1,1)$.  We show that of these four cases only the last
one is possible.  The \ul{first} case (see Figure~2a) cannot occur
because then $(0,1)\in S$\/ has for its $L-S$\/ neighbor $(1,1)$ and
a second point not adjacent to either $(1,1)$ or $x$, which must
thus be $y=(-1,2)$; similarly the second $L-S$\/ neighbor of $(-1,0)$
must be $y'=(-2,1)$; but then $L-S$\/ contains the adjacent $y,y'$,
contradiction.  The \ul{second} case (Fig.~2b) we dispose of similarly:
$(\pm1,0)$ are also the two $L-S$\/ neighbors of $x'=(0,1)$, so the
remaining six neighbors of~$x'$, including $x''=(1,1)$, are in~$S$\/;
thus $L-S$ contains a point adjacent to $x''$ but not to $x'$ or
$(1,0)$, which must be $(2,2)$.  Likewise $(2,-2)\in L-S$.  But
then $x_1,x_2=(2,\pm1)\in S$\/ each already has two $L-S$\/ neighbors,
namely $(1,0)$ and $(2,\pm2)$, so all of their other neighbors,
including $(2,0)$, $(3,0)$, and $(3,\pm1)$, are in~$S$.  But then
$(2,0)\in S$\/ (marked ``!''\ in Fig.~2b)
has only one neighbor in $L-S$, which is impossible.
The \ul{third} case reduces to the second, since $(\pm1,1)$ are then
also the two $L-S$\/ neighbors of $(0,1)$.  Thus the \ul{fourth} case
must obtain.  Then (see Fig.~2c) $L-S$\/ must also contain $(-2,-1)$
to account for $x_1=(-1,0)\in S$, and then also $(0,-2)$ to take care
of $x_2=(-1,-1)$, and so forth; inductively we find that $L-S\ni y \iff
y_1+2y_2\equiv1\bmod 5$.  Applying an arbitrary isometry of~$L$\/
to the resulting~$S$\/ recovers the ten~$S$\/ described in the
statement of the Corollary.~~\Qed

$$
\sq{
\1{}\1y\b\b\\ \1{y'\!}\b\b\0\\ \b\b\x\b\\ \b\0\b\b\\
\\ \\ \mbox{Figure 2a}
}
\qquad
\sq{
\b\b\b\0\b\\ \b\1{x'\!}\1{\!x''\!\!}\1{x_1\!}\b\\ \0\x\0\1!\b\\
\b\b\b\1{x_2\!}\b\\ \b\b\b\0\b\\ \\ \mbox{Figure 2b}
}
\qquad
\sq{
\b\0\b\b\\ \b\1{x_1\!}\x\0\\ \0\1{x_2\!}\b\b\\ \b\b\0\b\\
\\ \\ \mbox{Figure 2c}
}
$$

We can also easily prove that some of the bounds of Lemma~1 for
small~$n$ are attained.  Thus for $n=0$ we claim that any co-clique
$S\subset L$ has density at most~$1/4$; but this is clear because we
can partition $L$ into $2\times 2$ squares, each of which may contain
at most one element of~$S$.  Note that this also means that the
inequality of Prop.~1 cannot be sharp for $n<6$: if $S\subset L$
has maximal degree~$n$ and density $8/(16-n)$ then $N_{L-S}(y)=0$
for almost all $y\in L-S$, but then $L-S$ has density at most
$1/4 < 3/11 \leqs 1 - 8/(16-n)$.  Nevertheless we can use the
method of Prop.~1 to prove that $\delta(2)=1/2$.  This is because
if some $N_S(y)\geqs7$ for some $y\in L-S$\/ then $y$ has at least one
neighbor $x\in S$\/ with $N_S(x)\geqs4$.  Thus if $d(S)<4$ then
$N_S(y)\leqs 6$ for each $y\in L-S$.  In particular if $d(S)\leqs2$
then the number of adjacent pairs $(x,y)$ in $\B_r$ with $x\in S$,
$y\notin S$ is both $\geqs(6-O(1/r))|\B_r\cap S|$ and
$\leqs 6 |\B_r - S|$, so $(2-O(1/r)) |\B_r\cap S| \leqs |\B_r|$
and $S$\/ has density $\leqs 1/2$.

But this still leaves the matter of $\delta(1)$.  It turns out that
the bound $\delta(1)\geqs 1/3$ of Lemma~1 is again sharp, but it is
not readily proved by any of our methods thus far.  It does, however,
follow easily from a Voronoi decomposition of $\R^2$ relative
to~$S$; in fact this approach lets us treat the three cases
$n=0,1,2$ uniformly and quickly, and with some more work will
also prove that $\delta(3)=1/2$ and thus prove Conjecture~A and
the still-Life conjecture.

Given any discrete nonempty subset $S\subset\R^2$, the usual
{\em Voronoi decomposition} of~$\R^2$ relative to~$S$\/ (see Figure~3a)
is $\R^2 = \cup_{x\in S} V_x$, where the {\em Voronoi cell} $V_x$ is
$$
V_x := \{ z \in \R^2: |z-x'| \geqs |z-x| {\rm\ for\ all}\ x'\in S \}.
$$
The union of these cells is all of $\R^2$ because
$z\in V_x$ where $x$ is one of the elements of~$S$\/ closest to~$z$.
Each $V_x$ is a closed convex subset of $\R^2$ whose interior
contains~$x$, and the $V_x$ intersect only on their boundaries,
with $V_x\cap V_{x'}$ contained in the ``perpendicular bisector''
$\beta_{x,x'} := \{ z \in \R^2 : |z-x| = |z-x'| \}$.
Note that $V_x$ can be defined equivalently as
$$
V_x := \bigcap_{x'\in S \atop x'\neq x} H_x(x'),
$$
where $H_x(x')$ is the closed half-plane
containing~$x$ determined by $\beta_{x,x'}$.
If the $V_x$ are of bounded of diameter, the fact
that the $V_x$ tile the plane means that if $S$\/ has a density
$\delta$ --- defined as before, with the cardinality of $\B_r$
replaced by its area $4r^2$ --- then $1/\delta$ is the average
over $x\in S$ (i.e.\ the limit as $r\ra\infty$ of the average
over $x\in \B_r \cap S$\/) of the area~$A_x$ of~$V_x$.

\vspace*{1ex}

\centerline{
\setlength{\unitlength}{10pt}
\begin{picture}(10,10)
\put(3,2){\makebox(0,0)[c]{$\bullet$}}
\put(7,2){\makebox(0,0)[c]{$\bullet$}}
\put(1,6){\makebox(0,0)[c]{$\bullet$}}
\put(9,6){\makebox(0,0)[c]{$\bullet$}}
\put(5,8){\makebox(0,0)[c]{$\bullet$}}
\put(5,1){\line(0,1){3.667}}
\put(4,5){\line(3,-1){1}}
\put(4,5){\line(-2,-1){4}}
\put(4,5){\line(-1,2){2}}
\put(6,5){\line(-3,-1){1}}
\put(6,5){\line(2,-1){4}}
\put(6,5){\line(1,2){2}}
\put(.5,-.8){\makebox(9,1)[c]{Figure 3a}}
\end{picture}
\qquad\quad
\begin{picture}(11,9)
\put(3,2){\makebox(0,0)[c]{$\bullet$}}
\put(7,2){\makebox(0,0)[c]{$\bullet$}}
\put(1,6){\makebox(0,0)[c]{$\bullet$}}
\put(9,6){\makebox(0,0)[c]{$\bullet$}}
\put(5,8){\makebox(0,0)[c]{$\bullet$}}
\put(5,1){\line(0,1){3}}
\put(-.5,2.5){\line(1,1){1.5}}
\put(1,4){\line(1,0){2}}
\put(3,4){\line(1,1){1}}
\put(5,4){\line(-1,1){2}}
\put(3,6){\line(0,1){2}}
\put(3,8){\line(-1,1){1.5}}
\put(10.5,2.5){\line(-1,1){1.5}}
\put(9,4){\line(-1,0){2}}
\put(7,4){\line(-1,1){1}}
\put(5,4){\line(1,1){2}}
\put(7,6){\line(0,1){2}}
\put(7,8){\line(1,1){1.5}}
\put(.5,-.8){\makebox(9,1)[c]{Figure 3b}}
\end{picture}
\qquad\quad
\begin{picture}(11,9)
\put(2,6){\makebox(0,0)[c]{$\bullet$}}
\put(8,6){\makebox(0,0)[c]{$\bullet$}}
\put(5,3){\line(0,1){6}}
\put(5,3){\line(1,-1){2}}
\put(5,3){\line(-1,-1){2}}
\put(5,9){\line(1,1){2}}
\put(5,9){\line(-1,1){2}}
 \put(4.7,2.7){\line(1,0){.6}}
 \put(4.4,2.4){\line(1,0){1.2}}
 \put(4.1,2.1){\line(1,0){1.8}}
 \put(3.8,1.8){\line(1,0){2.4}}
 \put(3.5,1.5){\line(1,0){3.0}}
 \put(3.2,1.2){\line(1,0){3.6}}
 \put(4.7,9.3){\line(1,0){.6}}
 \put(4.4,9.6){\line(1,0){1.2}}
 \put(4.1,9.9){\line(1,0){1.8}}
 \put(3.8,10.2){\line(1,0){2.4}}
 \put(3.5,10.5){\line(1,0){3.0}}
 \put(3.2,10.8){\line(1,0){3.6}}
\put(.5,-.8){\makebox(9,1)[c]{Figure 3c}}
\end{picture}
}

\vspace*{2ex}

We called this the ``usual'' Voronoi decomposition because it uses
the usual Euclidean norm on $\R^2$.  But we could have used any
norm $|\cdot|$ to define the $V_x$, and for our purposes we would
like to use the $l_\infty$ norm: $|v|_\infty = \max(|v_1|,|v_2|)$
(See Fig.~3b).  The problem with that norm is that the associated
Voronoi cells can have non-boundary intersections, because if $x,x'$
have a coordinate in common then $\beta_{x,x'}$ has positive, indeed
infinite, area (see Fig.~3c).  Since this will invalidate the
description of $1/\delta$ as the area of the average Voronoi cell,
we must modify our definition.  Fortunately there is a natural
way to do this: when $x,x'$ have a coordinate in common, replace
$\beta_{x,x'}$ by the Euclidean perpendicular bisector of $x,x'$,
and use these $\beta_{x,x'}$ to define the $H_x(x')$ in the formula
$V_x = \cap_{x'} H_x(x')$.  Equivalently, if $x_1<x'_1\leqs x_2$
we regard the point $(\pm x_1,\pm x_2)$ as infinitesimally closer
to the origin than $(\pm x'_1,\pm x_2)$ by modifying the $l_\infty$
norm to a ``norm'' taking values in the nonnegative linear combinations
of~$1$ and the positive infinitesimal~$\eps$:
$$
|v|_\eps := \max(|v_1|,|v_2|) + \eps \min(|v_1|,|v_2|).
$$
We then define $V_x$ as the closure of
$$
\{ z \in \R^2: |z-x'|_\eps \geqs |z-x|_\eps {\rm\ for\ all}\ x'\in S \}.
$$
Either way, we find that $V_x$ is then a convex region containing~$x$
bounded by lines either parallel either to a coordinate axis or to
a line $x_2=\pm x_1$, and that these $V_x$ have disjoint interiors
and cover $\R^2$.

Moreover, if $S$\/ is a lattice subset then our apparently continuous
Voronoi construction actually yields a discrete object:
the boundary lines of the $V_x$ are all of the form $2x_i=c$ or
$x_1\pm x_2=c$ for some $c\in\Z$, so the $V_x$ are unions of isosceles
triangles from a fixed tiling of the plane (Fig.~4a).  These triangles
are the images under the isometry group Aut$(L)$ of the triangle
$\{ x: 0\leqs x_2\leqs x_1\leqs 1/2 \}$ of area~$1/8$.
[Indeed this triangle is a fundamental domain for Aut$(L)$.]
Thus the area $A_x$ of $V_x$ is a multiple of $1/8$ for any~$x$
for which $V_x$ is bounded.  But for our purposes we may assume
that the $V_x$ are of uniformly bounded diameter:
partition $L$ into $3\times 3$ squares~$Q$, and for any square
such that $Q\cap S = \emptyset$ augment $S$\/ by the center of~$Q$.
This does not decrease the upper density of~$S$, and does not
increase its maximum degree since all the added points are isolated;
and now every point of~$\R^2$ is at $l_\infty$ distance at most~$3$
from some point of~$S$.

\vspace*{1ex}

\centerline{
\setlength{\unitlength}{10pt}
\begin{picture}(10,11)
\multiput(1,2)(2,0){5}{ \multiput(0,0)(0,2){5}{ \circle*{.4} }}
\multiput(0,2)(2,0){5}{ \multiput(0,0)(0,2){5}{ \line(1,0){2} }}
\multiput(1,1)(2,0){5}{ \multiput(0,0)(0,2){5}{ \line(0,1){2} }}
\multiput(0,1)(2,0){5}{ \multiput(0,0)(0,2){5}{ \line(1,1){2} }}
\multiput(0,3)(2,0){5}{ \multiput(0,0)(0,2){5}{ \line(1,-1){2} }}
\multiput(2,1)(2,0){4}{ \line(0,1){10} }
\multiput(0,3)(0,2){4}{ \line(1,0){10} }
\put(.5,-.8){\makebox(9,1)[c]{Figure 4a}}
\end{picture}
\qquad\qquad\qquad
\begin{picture}(10,11)
\thicklines
\put(0,8){\circle{.4}} \put(2,8){\circle{.4}} \put(4,8){\circle{.4}}
\put(0,6){\circle{.4}} \put(2,6){\circle*{.4}} \put(4,6){\circle{.4}}
\put(0,4){\circle{.4}} \put(2,4){\circle{.4}} \put(4,4){\circle{.4}}
\put(2,5.3){\makebox(0,0)[c]{$x$}}
\put(0,4){\framebox(4,4){}}
\put(7,8){\circle*{.4}} \put(9,8){\circle{.4}} \put(11,8){\circle{.4}}
\put(7,6){\circle{.4}} \put(9,6){\circle*{.4}} \put(11,6){\circle*{.4}}
\put(7,4){\circle{.4}} \put(9,4){\circle{.4}} \put(11,4){\circle{.4}}
\put(9,5.3){\makebox(0,0)[c]{$x$}}
\put(7,4){\line(1,0){3}}
\put(10,4){\line(0,1){4}}
\put(10,8){\line(-1,0){1}}
\put(9,8){\line(-1,-1){2}}
\put(7,6){\line(0,-1){2}}
\thinlines
\put(7,4){\framebox(4,4){}}
 \put(10.25,4){\line(0,1){4}}
 \put(10.50,4){\line(0,1){4}}
 \put(10.75,4){\line(0,1){4}}
 \put(7,6.5){\line(1,0){0.5}}
 \put(7,6.8){\line(1,0){0.8}}
 \put(7,7.1){\line(1,0){1.1}}
 \put(7,7.4){\line(1,0){1.4}}
 \put(7,7.7){\line(1,0){1.7}}
\put(.5,-.8){\makebox(10,1)[c]{Figure 4b}}
\end{picture}
}

\vspace*{2ex}

We may thus bound $\delta$ from above by finding a lower bound on
the average Voronoi-cell area.  If $x$ is an isolated point of~$S$\/
then $V_x$ contains the $2\times 2$ square
$\{ z: |z-x|_\infty \leqs 1 \}$, so $A_x\geqs 4$; therefore if
$S$\/ consists entirely of isolated points its density is at most
$1/4$, and we have recovered the bound $\delta(0)=1/4$.  More generally,
for each $x\in S$ we have $A_x \geqs 4 - N_S(x)$, because each
neighbor of~$x$ cuts out from that $2\times 2$ square either
a $2\times (1/2)$ rectangle of area~$1$ (if it is orthogonally
adjacent) or an isosceles right triangle of area~$1/2$ (if diagonally
adjacent; see Fig.~4b).  We deduce:

{\bf Proposition 2.}
{\sl For each $n$ we have $\delta(n)\leqs 1/(4-n)$.
Thus (by Lemma~1) we have $\delta(n) = 1/(4-n)$ for $n=0,1,2$.
Moreover a periodic set $S$ with $d(S)\leqs n$ has density
$1/(4-n)$ if and only if one of the following holds:
\begin{description}
\romit{\qquad i)} $n=0$, and the $2\times 2$ squares centered at
    points of~$S$ tile the plane;
\romit{\qquad ii)} $n=1$, and $S$ consists of dominos, i.e.\ of pairs
    $\{ x, x+(1,0) \}$ or $\{ x, x+(0,1) \}$, whose circumscribed
    $3\times2$ or $2\times3$ rectangles tile the plane;
\romit{\qquad iii)} $n=2$, and $S$ is one of the four equivalent sets
    $\{ x: x_1\equiv0\bmod2 \}$, $\{ x: x_1\equiv1\bmod2 \}$,
    $\{ x: x_2\equiv0\bmod2 \}$, $\{ x: x_2\equiv1\bmod2 \}$.
\end{description}
}

{\em Proof}\/: For $x\in S$ let $V^0_x$ be the intersection of $V_x$
with the $2\times 2$ square $\{ z: |z-x|_\infty \leqs 1 \}$ centered
at~$x$ --- we shall call this the {\em narrow Voronoi cell} about~$x$
--- and let $A^0_x\leqs A_x$ be the area of $V^0_x$.
Then $V^0_x$, and thus also $A^0_x$, depends only on the neighbors
of~$x$ in~$S$.  We have already seen that $A^0_x\geqs 4-N_S(x)$, from
which $\delta(n)\leqs1/(4-n)$ follows.  We readily check that $A^0_x$
is as small as $4-N_S(x)$ only in three cases:
(i) when $x$ is isolated; (ii) when $x$ has one orthogonally adjacent
neighbor in~$S$\/; or (iii) when $x$ has two orthogonally adjacent
neighbors on opposite sides.  Thus a periodic set $S$\/ of maximal
degree $n$ has density $1/(4-n)$ if and only if $n=0$, 1, or~2,
with each $x\in S$\/ having $n$ neighbors distributed as indicated
in the previous sentence under (i), (ii), or (iii) respectively,
and the $V^0_x$ tiling $\R^2$.  For $n=0$ this is exactly case (i)
of the Proposition.  For $n=1$, note that the narrow Voronoi cells
of the two elements of a domino coalesce to the $3\times 2$ or
$2\times 3$ rectangle inscribing it.  Finally for $n=2$ the condition
that each $x$ have two orthogonally adjacent and opposite neighbors
forces all the $x\in S$\/ to line up in rows or columns, and it
is then clear that $S$\/ can only attain density $1/2$ if these
rows and columns are spaced two units apart.~~\Qed

As in the Corollary to Prop.~1 we may also describe the periodic sets
achieving $\delta(0)=1/4$.  These are
$\{ x \in L: x_1 \equiv c \bmod 2, x_2 \equiv b(x_1) \bmod 2 \}$
and $\{ x \in L: x_2 \equiv c \bmod 2, x_1 \equiv b(x_2) \bmod 2 \}$
for some integer $b \bmod 3$ and some periodic function
$b: 2\Z+c \ra \Z/2$; that is, we start with the obvious tiling
of the plane by $2\times 2$ squares and then periodically shift
some columns or rows by one unit.  The description of periodic
sets achieving $\delta(1)=1/3$ is more complicated.
Besides shifting the obvious tilings with $3\times2$ or
$2\times3$ rectangles we may use both $3\times 2$ and $2\times 3$
layers in the same pattern as long as the choices of orientation
as well as shift are periodic; there are also many periodic tilings
which do not decompose into horizontal or vertical layers at
all.\footnote{Thank to Dan Hoey for pointing this out; some of these
tilings are loosely analogous to the ``tv static'' illustrated
in Appendix~A.}

From our results thus far it follows too that $\delta(2),\delta(3)$
are the only possible exception to our expectation that $\delta(n)$
be strictly increasing:

{\bf Corollary.} {\sl For each positive $n\leqs8$ with the possible
exception of $n=3$ we have $\delta(n)>\delta(n-1)$.}

{\em Proof}\/:  For $n=1,2,7,8$ this is clear from the values
of $\delta(n)$ already known.  So it remains to prove
$\delta(3)<\delta(4)<\delta(5)<\delta(6)$.  We estimate
$\delta(3)$ above as we did for $\delta(2)$: if $d(S)\leqs 3$
then $N_S(y)\leqs 6$ for all $y\notin S$\/ so $S$ has density
at most $6/(6+5)=6/11$.  This combined with the inequalities
of Lemma~1 and Prop.~1 yields
$$
\delta(3) \leqs \frac6{11} < \frac35 \leqs \delta(4) \leqs \frac23
< \frac9{13} \leqs \delta(5) \leqs \frac8{11} < \frac45 = \delta(6),
$$
and we are done.~~\Qed

{\bf 3. Digression: bounds on $\delta(4)$ and $\delta(5)$}

The upper bounds $2/3$, $8/11$ on $\delta(4)$, $\delta(5)$ are not
the best known.  For several years the records were held by
Greg Kuperberg~\cite{Greg}; I describe them here with his permission.
The method is as follows: for some finite region $B\subset L$,
find $D>0$ and nonnegative weights $w_x$ ($x\in B$) such that
any $S\in L$ with $d(S)\leqs n$ satisfies the inequality
$\sum_{x\in S} w_x \leqs D \sum_{x\in B} w_x$.  Then
by averaging over translates an arbitrary set~$S$\/ of maximum
degree $\leqs n$ we find that $S$\/ has density at most~$D$\/;
therefore $\delta(n)\leqs D$.  Given $B$, an optimal system
of weights $w_x$ may be found by linear programming.  This generalizes
one of the methods we used above; for instance the first argument for
$\delta(0)\leqs1/4$ amounts to taking
\hbox{$B=\{x: 0\leqs x_i\leqs 1\}$,}
$w_x=1$ for all $x\in B$, and $D=1/4$.  But one can do better
with varying weights; for instance, G.~Kuperberg computed that the
following weights on a $6\times 6$ square yield $\delta(4)\leqs 5/8$:
$$
\begin{array}{cccccc}
0 & 1 &   3  &   3  & 1 & 0 \\
1 & 5 &   8  &   8  & 5 & 1 \\
3 & 8 &\!13\!&\!13\!& 8 & 3 \\
3 & 8 &\!13\!&\!13\!& 8 & 3 \\
1 & 5 &   8  &   8  & 5 & 1 \\
0 & 1 &   3  &   3  & 1 & 0
\end{array}
\ \ .
$$
Appealing as the Fibonacci coefficients and simple fraction $5/8$
are, they are not optimal: Kuperberg proved (by programming a computer
to list all patterns attaining
$\sum_{x\in S} w_x \leqs (5/8) \sum_{x\in B} w_x$) that
$\delta(4)$ must be strictly smaller than $5/8$, and searched
unsuccessfully for periodic sets of any density $>3/5$.  As a
byproduct he did find that there are many more ways of attaining
$3/5$ by locally perturbing the ``alternative pattern'' of Appendix~B
(see Figure~5a); a particularly nice pattern is obtained by applying
this perturbation as densely as possible (Fig.~5b).
$$
\renewcommand{\arraystretch}{.9}
\def\b{\makebox[.36cm]{$\bullet$}}
\def\:{\makebox[.36cm]{$\cdot$}}
\sq{
\:\b\:\b\:\makebox[1.5cm]{}\:\b\b\b\:\\
\b\b\:\b\b\makebox[1.5cm]{$\longrightarrow$}\b\:\:\:\b\\
\:\b\:\b\:\makebox[1.5cm]{}\:\b\b\b\:\\
\\ \mbox{Figure 5a}
}
$$
\vspace*{2ex}
$$
\renewcommand{\arraystretch}{.8}
\def\b{\makebox[.32cm]{$\bullet$}}
\def\:{\makebox[.32cm]{$\cdot$}}
\sq{
\:\:\b\b\b\:\:\:\b\b\b\:\:\:\b\b\b\:\:\:\\
\b\b\b\:\b\b\b\b\b\:\b\b\b\b\b\:\b\b\b\b\\
\:\:\b\:\b\:\:\:\b\:\b\:\:\:\b\:\b\:\:\:\\
\b\b\b\:\b\b\b\b\b\:\b\b\b\b\b\:\b\b\b\b\\
\:\:\b\b\b\:\:\:\b\b\b\:\:\:\b\b\b\:\:\:\\
\b\b\:\:\:\b\b\b\:\:\:\b\b\b\:\:\:\b\b\b\\
\:\b\b\b\b\b\:\b\b\b\b\b\:\b\b\b\b\b\:\b\\
\:\b\:\:\:\b\:\b\:\:\:\b\:\b\:\:\:\b\:\b\\
\:\b\b\b\b\b\:\b\b\b\b\b\:\b\b\b\b\b\:\b\\
\b\b\:\:\:\b\b\b\:\:\:\b\b\b\:\:\:\b\b\b\\
\:\:\b\b\b\:\:\:\b\b\b\:\:\:\b\b\b\:\:\:\\
\b\b\b\:\b\b\b\b\b\:\b\b\b\b\b\:\b\b\b\b\\
\:\:\b\:\b\:\:\:\b\:\b\:\:\:\b\:\b\:\:\:\\
\b\b\b\:\b\b\b\b\b\:\b\b\b\b\b\:\b\b\b\b\\
\:\:\b\b\b\:\:\:\b\b\b\:\:\:\b\b\b\:\:\:\\
\b\b\:\:\:\b\b\b\:\:\:\b\b\b\:\:\:\b\b\b\\
\:\b\b\b\b\b\:\b\b\b\b\b\:\b\b\b\b\b\:\b\\
\:\b\:\:\:\b\:\b\:\:\:\b\:\b\:\:\:\b\:\b\\
\:\b\b\b\b\b\:\b\b\b\b\b\:\b\b\b\b\b\:\b\\
\b\b\:\:\:\b\b\b\:\:\:\b\b\b\:\:\:\b\b\b\\
\\ \mbox{Figure 5b}
}
\renewcommand{\arraystretch}{1}
\def\b{\makebox[.4cm]{$\bullet$}}
\def\:{\makebox[.4cm]{$\cdot$}}
$$
\vspace*{-20ex}

For $\delta(5)$ the upper bound $149/212$ is obtained from the weights
$$
\begin{array}{cccccc}
1 & 2 & 3 & 3 & 2 & 1\\
2 & 6 &10 &10 & 6 & 2\\
3 &10 &16 &16 &10 & 3\\
3 &10 &16 &16 &10 & 3\\
2 & 6 &10 &10 & 6 & 2\\
1 & 2 & 3 & 3 & 2 & 1
\end{array}
\ \ ;
$$
numerically $149/212=.70283+$ while $9/13=.69231-$.  As $N$\/
increases, the bounds obtained from an $N\times N$\/ square must
decrease to $\delta(n)$, but of course they become exponentially
more difficult to compute (because each maximal admissible subset
of the square puts a linear relation on the $w_x$), though $N=7$
and maybe even $N=8$ should be accessible to current technology.
Unfortunately for $n=4$ and $n=5$ Kuperberg shows that if the
conjectured values $3/5$ and $9/13$ on $\delta(n)$ are correct
then no finite $N$\/ will suffice to prove them.

More recently (Dec.~95) Marcus Moore announced on the {\sc usenet}
newsgroup {\tt comp.theory.cell-automata} a new approach to
bounding $\delta(n)$ from above, closely related to ``discharging''
arguments as explained in~\cite{SK} (thanks to Allan Wechsler
for this reference).  In that forum, and in later personal
communication~\cite{Moore}, Moore claims a simpler proof of
$\delta(3)=1/2$, as well as the new results $\delta(5)=9/13$
and $\delta(4)\leqs 8/13$.  To my knowledge these results
have yet to appear even in preprint form.

\vspace*{3ex}

{\bf 4.~Proof of Conjecture A.}

If $d(S)\leqs3$ we can no longer assert that every Voronoi cell
has area at least~2.  For instance if $S$\/ contains $x=(0,0)$
and its three neighbors $(0,1)$ and $(\pm1,0)$ then $A^0_x=3/2$
and $A_x$ need be no larger.  Note, however, that in this case
the Voronoi cell of $(0,1)\in S$\/ must have area at least $11/4$.
In general we shall see that S can be partitioned into uniformly
bounded subsets $S_0$ on which $A_x$ averages to at least 2.
(This method is thus also related, though more loosely than Moore's,
with the ``discharging'' techniques of~\cite{SK}: it can be regarded
as reapportioning the $A_x$ locally to normalized cell areas
$A'_x\in\frac18\Z$ which still average to the inverse density of~$S$\/
but satisfy $A'_x\geqs2$ for all~$x$.  The method is even more strongly
reminiscent of Hsiang's approach to Kepler's conjecture, though
fortunately not as many cases need be considered\ldots)\ \
Thus we will settle Conjecture~A by proving

{\bf Theorem.}  {\sl $\delta(3)=1/2$.  Moreover every periodic set
$S\subset L$ of density $1/2$ with $d(S)\leqs3$ has $N_S(x)\geqs2$
and $N_S(y)\geqs4$ for all $x\in S$, $y\notin S$, so in particular
$S$ is automatically a periodic still Life of maximal density.}

{\em Proof}\/: As before we may assume that all the Voronoi cells
$V_x$ ($x\in S$\/) have finite area $A_x$.  For $x\in S$\/ let 
$\alpha(x),\alpha^0(x)$ be the integers $8(A_x-2)$ and
$8(A^0_x-2)\leqs\alpha(x)$.  Note that $\alpha^0$ depends
only on the neighbors of~$x$ in~$S$.  In the course of proving Prop.~2
we noted in effect that if $N_S(x)\leqs1$ then $\alpha^0(x)=8$,
12, or~16 according as $x$ has one orthogonal neighbor, one
diagonal neighbor, or no neighbors.  For $N_S(x)=2$ and $N_S(x)=3$
there are respectively six and ten configurations up to symmetry;
we label them and list their $\alpha^0$ values (using $\bullet$ and
$\circ$ for elements of~$S$\/ and $L-S$\/ respectively):

$$
L: \ \sq{\0\0\0\\ \b\x\b \\ \0\0\0} \ \alpha^0(x) = 0
\qquad\qquad
L': \ \sq{\0\b\0\\ \0\x\b \\ \0\0\0} \ \alpha^0(x) = 2
$$ $$
M: \ \sq{\0\0\b \\ \b\x\0 \\ \0\0\0 } \ \alpha^0(x) = 4
\qquad\qquad
M': \ \sq{\0\0\b \\ \0\x\b \\ \0\0\0} \ \alpha^0(x) = 7
$$ $$
N: \ \sq{\b\0\0 \\ \0\x\0 \\ \0\0\b } \ \alpha^0(x) = 8
\qquad\qquad
N': \ \sq{\0\0\b \\ \0\x\0 \\ \0\0\b} \ \alpha^0(x) = 8
$$
\vspace*{1ex}
$$
A: \ \sq{\0\b\0 \\ \b\x\b \\ \0\0\0 } \ \alpha^0(x) = -4
\qquad\qquad
D: \ \sq{\b\0\b \\ \0\x\0 \\ \0\0\b} \ \alpha^0(x) = 4
$$ $$
B: \ \sq{\0\0\b \\ \b\x\b \\ \0\0\0 } \ \alpha^0(x) = -1
\qquad\qquad
E: \ \sq{\b\b\0 \\ \0\x\0 \\ \0\0\b} \ \alpha^0(x) = 3
$$ $$
C: \ \sq{\0\b\0 \\ \b\x\0 \\ \0\0\b } \ \alpha^0(x) = -2
\qquad\qquad
F: \ \sq{\b\0\b \\ \0\x\0 \\ \0\b\0} \ \alpha^0(x) = 0
$$ $$
C': \ \sq{\b\b\0 \\ \0\x\b \\ \0\0\0 } \ \alpha^0(x) = 1
\qquad\qquad
F': \ \sq{\b\0\b \\ \0\x\b \\ \0\0\0} \ \alpha^0(x) = 3
$$ $$
C'': \ \sq{\0\b\b \\ \0\x\b \\ \0\0\0 } \ \alpha^0(x) = 2
\qquad\qquad
F'': \ \sq{\b\b\b \\ \0\x\0 \\ \0\0\0} \ \alpha^0(x) = 6
$$
We shall show that $S$\/ can be partitioned into
subsets $S_0$ on which $\alpha$ has nonnegative average,
and which are uniformly bounded, i.e.\ are each contained in
a translate of a fixed bounded set in the plane ($\B_{10}$
is more than large enough).
For a finite subset $S_0\subset S$\/ we will denote by $\sigma S_0$
and $\sigma^0 S_0$ the sums
$$
\sigma S_0 := \sum_{x\in S_0} \alpha(x),
\qquad
\sigma^0 S_0 := \sum_{x\in S_0} \alpha^0(x) \leqs \sigma S_0.
$$

First, each $A$ element has a unique $F''$ neighbor and vice versa;
we pair these, and note that their $\alpha^0$'s add to $6-4=2$.
This leaves only the $B$\/ and $C$\/ elements.  For each $x\in S$\/
of type~$B$\/ or~$C$\/ there is a unique pair $x',x''$ of neighbors
of~$x$\/ in~$S$\/ adjacent to each other; we group $x$ with $x'$
and $x''$.  Note that a $C$\/ may be adjacent to one or two $B$\/'s,
in which case they are combined for two or three reasons;
and that a $B$\/ may be joined with an existing $AF''$ pair.
Also, an $AF''$ pair, or two orthogonally adjacent $C'$ elements,
may be combined with $B$\/'s on both sides, in which case
we consider those $B$\/'s and their two
common neighbors as a four-element grouping (``quartet'')
with a $\sigma^0$ of~$0$:
$$
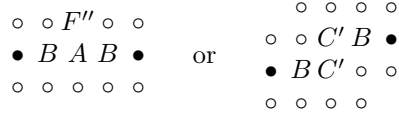
\sq{\0\0\1{F''}\0\0 \\ \b\1B\1A\1B\b \\ \0\0\0\0\0}
\quad {\rm or} \quad
\sq{\1{}\0\0\0\0 \\ \0\0\1{C'}\1B\b \\ \b\1B\1{C'}\0\0 \\ \0\0\0\0\1{}}
$$

If $N_S(x')\leqs2$ then $x'$ has type $L'$ or $M'$.  If $x'$ is of
type~$L'$ then $x$ is of type $B$\/ and $x''$ is of type $B$, $E$, $F'$,
or~$M$\/; in the first case, $\sigma^0\{ x, x', x'' \} = 0$:
$$
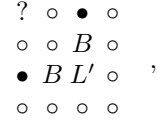
\sq{ \1?\0\b\0 \\ \0\0\1B\0 \\ \b\1B\1{L'}\0 \\ \0\0\0\0 } \ ,
$$
 (recall that ``?''\ means a point which may be either in $S$\/
 or $L-S$\/),
and otherwise, $\sigma^0\{ x, x', x'' \} \geqs 4$.
If $x'$ is of type $M'$ then again
$\sigma^0\{ x, x', x'' \} \geqs 4$, with equality iff
$\{ x, x', x'' \}$ is a $BCM'$ grouping.

While the $BL'B$ grouping already has narrow Voronoi cells of
average area~$2$, we later need to use the large Voronoi cells
of an $L'$ point elsewhere.  Fortunately we can remove the $L'$
point $x'$ from the grouping and compensate for the negative
$\alpha^0$'s of the $B$\/ points $x,x''$ by analyzing the point~$z$
marked ``?''\ in the $BL'B$\/ diagram.  If $z\notin S$\/
then the Voronoi cells of $x,x''$ extend beyond $V^0_x,V^0_{x''}$
into the square of side~1 centered at~$z$; thus
$\alpha(x),\alpha''(x)\geqs0$.  Otherwise $z\in S$\/ and
$\alpha(z)>0$.  Moreover $\alpha(z)\geqs2$ unless $z$ is of
type~$C'$, in which case one of the neighbors of~$z$ in~$S$\/
has type~$F'$:
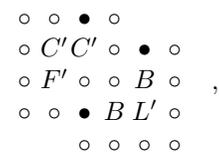
$$
\sq{
\0\0\b\0\1{}\1{}\\
\0\1{C'}\1{C'}\0\b\0 \\
\0\1{F'}\0\0\1B\0 \\
\0\0\b\1B\1{L'}\0 \\
\1{}\1{}\0\0\0\0
}
\ ,
$$
in which case we group $x,x''$ with that neighbor.  Otherwise we group
$x,x''$ with~$z$.  Then $\sigma^0 \{ x, x'', z \} \geqs -1 -1 +2 = 0$;
but we cannot yet conclude that our grouping has nonnegative $\sigma^0$,
because $z$ might be the ``?''-point of more than one $BL'B$
configuration, or may be already involved in another grouping, or both.
The preexisting groupings are $AF''$ pairs, $B$\/ and $C$\/ triangles,
$BB$\/ pairs obtained by deleting the $L'$ member of a $BL'B$ triangle,
and $BAF''B$, $BC'C'B$ quartets.
Clearly $z$ may not be contained in a quartet or a $BB$\/ pair.
Thus if $z$ is of type $F''$ then $z$ is part of an $AF''$ pair
or a $BAF''$ triangle.  Let $z',z''$ be the diagonal neighbors
of~$z$ in~$S$, with $z'$ closer to the $BL'B$ triangle.  Then
$z'$ is of type $F'$:
$$
\sq{
\0\1{z''}\0\0\1{}\1{}\\
\0\1{\!A\,}\1{F''}\0\b\0 \\
\0\1{F'}\0\0\1B\0 \\
\0\0\b\1B\1{L'}\0 \\
\1{}\1{}\0\0\0\0
}
\ ,
$$
with $\alpha^0(z')=3$.  Since $z'$ is not yet part of a grouping,
we attach it to $x,x''$ and our $AF''$ pair, and thus also to $z''$
if $z''$ is of type~$B$.  If $z$ is the ``?''-point of a second
$BL'B$ triangle then $z''$ is of type $F''$ and we have an octet
of four $B$\/'s, two $F'$'s, an $A$ and an~$F''$:
$$
\sq{
\1{}\1{}\0\0\0\0 \\
\0\0\b\1{\bf B}\1{L'}\0 \\
\0\1{\bf F'}\0\0\1{\bf B}\0 \\
\0\1{\bf A}\1{\bf F''}\0\b\0 \\
\0\1{\bf F'}\0\0\1{\bf B}\0 \\
\0\0\b\1{\bf B}\1{L'}\0 \\
\1{}\1{}\0\0\0\0
}
\ .
$$
 Otherwise we have a quintet or sextet of $A$, $F$, $F'$, and
 two or three $B$\/s.
Thus whenever $z$ is of type $F''$ we incorporate it into a grouping
with a positive $\sigma^0$.  If $z$ is not of this type, but
is the ``?''-point of two, three or four $BL'B$\/ triangles, then
$\alpha^0(z)$ is at least 7, 8, or 16 respectively, more than
enough to compensate for the $-4$, $-6$ or $-8$ of the $B$\/'s since
in this case $z$ is either previously unattached or, as in the
next diagram, part of a $CM'F'$ triangle with a $\sigma^0$ of~8:
$$
\sq{
\1{}\1{}\0\0\0\0 \\
\1{}\1{}\b\1B\1{L'}\0 \\
\b\0\0\0\1B\0 \\
\0\1C\1{M'}\0\b\0 \\
\0\1{F'}\0\0\1B\0 \\
\0\0\b\1B\1{L'}\0 \\
\1{}\1{}\0\0\0\0
}
\ .
$$
Assume then that $z$ is not of type $F''$, nor the ``?''-point of
more than one $BL'B$\/ triangle, but is nevertheless already involved
 in a grouping.  Such a grouping must be a $B$\/ or $C$\/~triangle,
 or a $BB$\/~pair augmented by an $AF'F''$ triangle and possibly
 an attached $B$\/ or $BBF'$ to a quintet, sextet or octet.
 In the latter case $z$ is of type $F'$ and becomes incorporated
 into an $AB^4F'F''$ septet with $\sigma^0=1$:
$$
\sq{
\1{}\1{}\0\0\b\0\1{}\1{} \\
\0\b\0\1{\bf F''}\1{\bf A}\b\1{}\1{} \\
\0\1{\bf B}\0\0\1{\bf F'}\0\b\0 \\
\0\1{L'}\1{\bf B}\1F\0\0\1{\bf B}\0 \\
\0\0\0\0\b\1{\bf B}\1{L'}\0 \\
\1{}\1{}\1{}\1{}\0\0\0\0
}
$$
or an $AB^5F'F''$ octet or $AB^8F'^2F''$ dozen with $\sigma^0=0$:
$$
\sq{
\1{}\1{}\1{}\0\b\0\1{}\1{} \\
\1{}\1{}\0\0\1{\bf B}\0\1{}\1{} \\
\0\b\0\1{\bf F''}\1{\bf A}\0\1{}\1{} \\
\0\1{\bf B}\0\0\1{\bf F'}\0\b\0 \\
\0\1{L'}\1{\bf B}\1F\0\0\1{\bf B}\0 \\
\0\0\0\0\b\1{\bf B}\1{L'}\0 \\
\1{}\1{}\1{}\1{}\0\0\0\0
}
\qquad\quad
\sq{
\1{}\1{}\1{}\1{}\0\0\0\0 \\
\0\0\0\0\b \1{\bf B}\1{L'}\0 \\
\0\1{L'}\1{\bf B}\1F\0\0\1{\bf B}\0 \\
\0\1{\bf B}\0\0\1{\bf F'}\0\b \0 \\
\0\1L\0\1{\bf F''}\1{\bf A}\0\1{}\1{} \\
\0\1{\bf B}\0\0\1{\bf F'}\0\b\0 \\
\0\1{L'}\1{\bf B}\1F\0\0\1{\bf B}\0 \\
\0\0\0\0\b \1{\bf B}\1{L'}\0 \\
\1{}\1{}\1{}\1{}\0\0\0\0
}
\;\ .
$$

 In the case of a $B$\/ or $C$\/ triangle, $z$ is of type $E$, $L'$
 or $M'$.  The $E$\/ case arises only for a $CEF'$ triangle, which
 we group with the $BB$\/~pair to raise its $\sigma^0$ to~$+2$.
 In the $L'$~case one of $z$'s
neighbors in~$S$\/ has type~$B$, and the other neighbor~$z'$ may
have type $E$, $F'$ or~$M'$ --- but may not be another $B$,
because in that case we have already deleted $z$ from the resulting
$BL'B$ triangle!  Thus $z$ is involved in a $B$\/~triangle with a
$\sigma^0$ of $7+2-1=8$ if $z'$ is of type $M'$ and
$3+2-1=4$ otherwise.  In the latter case we include the points $x,x''$
of type~$B$\/ to get a quintet with a $\sigma^0$ of~$2$.  In
the former case $z'$ may itself be a ``?''-point so we may include
either two or four points of type~$B$, but the resulting quintet
or septet still has a positive (indeed $\geqs4$) $\sigma^0$.
Finally if $z$ is of type $M'$ then it is part of a $BL'M'$, $BCM'$,
$ECM'$, $F'CM'$ or $M'CM'$ triangle two which we will incorporate two,
or possibly (in the first and last case) four, points of type~$B$,
and in this case two we have put $x,x'',z$ in a quintet or septet
with a $\sigma^0$ of at least~$2$ or~$4$ respectively.

Thus for each $BB$-pair either $z\notin S$, so the pair already
has nonnegative $\sigma^0$, or $z\in S$ and we have succeeded
in extending the pair by $z$ and perhaps other points of~$S$\/
to a subset $S_0\subset S$\/ contained in a translate of $\B_4$,
with $|S_0|\leqs12$ and $\sigma^0(S_0)\geqs0$.
We note for future reference that $\sigma^0(S_0)$ vanishes only when
 either $z$ is of type $C''$ or $L'$ and $S_0=\{x,x'',z\}$ or $z$
 is of type~$F'$ and contained in an $AB^5F'F''$ octet or $AB^8F'^2F''$
 dozen with $\sigma^0=0$, and that if
$S_0$ contains the diagonal neighbor of a type-$C$\/ element of~$S$\/
then $\sigma^0(S_0)\geqs2$ with equality only if $z$ is of type~$D$\/:
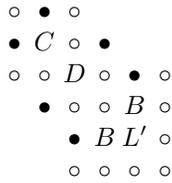
$$
\sq{
\0\b\0\1{}\1{}\1{} \\
\b\1C\0\b\1{}\1{} \\
\0\0\1D\0\b\0 \\
\1{}\b\0\0\1B\0 \\
\1{}\1{}\b\1B\1{L'}\0 \\
\1{}\1{}\0\0\0\0
}
$$

We are left with the cases in which $x$ is of type $B$ or~$C$\/
and $N_S(x')=N_S(x'')=3$ but $x,x',x''$
are not part of a $BAF''B$\/ or $BC'C'B$\/ quartet.  Of these,
several have $\sigma^0 \{ x, x', x'' \} > 0$, and two
configurations have $\sigma^0 \{ x, x', x'' \} = 0$:
$$
\sq{ \1{}\b\0\0 \\ \0\0\1{F'}\0 \\ \b\1B\1C\0 \\ \0\0\0\b }
\quad {\rm and} \quad
\sq{ \1{}\0\0\b \\ \0\0\1E\0 \\ \b\1B\1C\0 \\ \0\0\0\b }
\ .
$$
There remains the $BCB$\/ configuration, with a $\sigma^0$ of~$-4$:
$$
\sq{ \1?\0\b\0 \\ \0\0\1B\0 \\ \b\1B\1C\0 \\ \0\0\0\1{x_1} }
\ .
$$
Again we remove the $C$\/ point, and analyze the remaining $BB$\/ 
pair as we did for $BL'B$ triangles. 
Let $z$ be the point marked ``?''\ in this diagram.
From our $BL'B$\/ analysis we know that if $z\notin S$\/
then $z$ contributes at least~1 to the $\alpha$ of each $B$~point,
and if $z\in S$\/ then we may add at most ten points around the
$BCB$\/ triangle to increase its $\sigma^0$ by at least~$2$.
(This is why we worked so hard to excise the $L'$~point of a
$BL'B$\/~triangle: $z$ might be such a point, and if so we need
$\alpha^0(z)=+2$ to balance the $-2$ of such the $BB$\/ pair in
a $BCB$\/~triangle.)

We are left with the $C$\/ point and its $\alpha^0$ of~$-2$.
Consider its diagonal neighbor~$x_1$.  If $\alpha^0(x_1)\leqs0$
then $x_1$ is of type~$F$\/ or~$C$\/:
$$
\sq{
\1{}\0\b\0\1{} \\
\0\0\1B\0\1{z'} \\
\b\1B\1C\0\0 \\
\0\0\0\1F\b \\
\1{}\1{}\b\0\0
}
\quad {\rm or} \quad
\sq{
\1{}\0\b\0\1{} \\
\0\0\1B\0\1{z'} \\
\b\1B\1C\0\0 \\
\0\0\0\1C\b \\
\1{}\1{z''}\0\b\0
}
\ .
$$
In either case we gain another $+2$ from the point marked $z'$,
to which our $BL'B$ analysis applies again.  Note that if $x_1$
is of type~$C$\/ then it may be part of a $BCB$\/ triangle, but
in that case we use the point $z''$ to raise the $\alpha$ of
the $x_1$ grouping and thus avoid counting the same $+2$
contribution twice.
If $\alpha^0(x_1)$ is positive then it and its associated grouping
have a total $\alpha^0$ of at least~2 (note that $x_1$ cannot be of
type~$B$\/), so by combining the groupings of $x_1$ and the $BCB$\/
triangle we obtain a grouping with $\sigma\leqs0$, unless $x_1$
is of type~$D$\/ and attached to three such $C$\/'s.
(If $x_1$ is attached to more than one $BCB$\/
triangle then, since $\alpha^0(x_1)>0$, the type of $x_1$ must be
either $D$, $N$, or~$N'$, and then $\alpha^0(x_1)\geqs4$ can handle
two $C$\/'s with their $\alpha^0$ of~$-2$, but the case of three
neighbors requires an additional argument.)  We then have the
following situation, in which the points marked $u,u',v$ may
be in~$S$\/ or~not:
$$
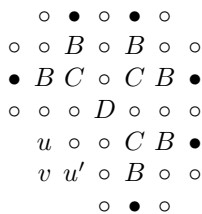
\sq{
\1{}\0\b\0\b\0\1{} \\
\0\0\1B\0\1B\0\0 \\
\b\1B\1C\0\1C\1B\b \\
\0\0\0\1D\0\0\0 \\
\1{}\1u\0\0\1C\1B\b \\
\1{}\1v\1{u'}\0\1B\0\0 \\
\1{}\1{}\1{}\0\b\0\1{}
}
$$
If $u\notin S$\/ then the Voronoi cells of $x_1$ and of its
adjacent $C$\/~point nearest to~$u$ extend by at least $1/4$ and $1/8$
into the square of side~1 centered at~$u$, increasing the $\sigma$ of
the $DCCC$ grouping around $x_1$ to at least $-2+3>0$. Likewise
we dispose of the case $u'\notin S$.
Thus we may assume that $u,u'\in S$.  We next ask
whether $v\in S$.  If not then $u,u'$ are of type $F'$ or worse
($M,$ $ M',$ $N$, or of degree $\leqs 1$).  All the groupings thus
far including such points have a $\sigma^0$ of at least twice the
number of available $F'$ points, with the exception of a $BCF'$ trio
withe $\sigma^0=0$; but $u,u'$ cannot both be involved in such trios
since the purported $B$\/ points would actually have type~$F'$.
Thus by grouping $u,u'$ with $x_1$ we obtain a grouping with
$\sigma^0\geqs0$.  If on the other hand $v\in S$\/ then at
least one of $u,u'$ has only two neighbors (else $N_S(v)\geqs 4$),
and thus has type $M'$; grouping that $M'$ point with $x_1$
finally eliminates the last possible grouping with negative $\sigma$.

This completes the proof of $\delta(3)=1/2$.  Moreover it characterizes
the periodic sets $S\subset L$ of density $1/2$ with $d(S)\leqs3$ as
the sets which when partitioned as above decompose into subsets $S_0$
with $\alpha(S_0)=0$.  By listing such subsets and their neighborhoods
we see that no element of~$S$\/ may have $N_S$ less than~1, and no
point of $L-S$\/ may have less than~4 neighbors in~$S$.  This completes
the demonstration of the Theorem.~~\Qed\,\Qed

Examining the neighborhoods of groupings with vanishing $\alpha(S_0)$
leads us to exclude several more possible $S_0$; we believe that
the still-Lifes of density~$1/2$ shown in Appendix~A show all
possible $S_0$.  In particular we assert that types $M,M',N,N',E$\/
cannot occur, and that a type-$D$\/ point must have one or two
neighbors of type~$C$\/ (the elegant pattern discovered by Moore
shows that first possibility, and the last pattern of Appendix~A shows
$DCC$\/ groupings).  The ``chicken wire''$\!$, ``tv~static''$\!$,
``square waves'' and ``Hoey~I'' patterns have many variations, and we
show them in quasiperiodic versions from which many periodic ones can
be extracted (cf.\ the descriptions of sets attaining $\delta(n)$
for $n=0,1,7$ in~\S2).  Moreover ``onion bulbs'' of order~$k$
may be interleaved with ``chicken wire'' of period~$4k$ or
(for $k=2$) with layers of ``Hoey~II'' and ``$DCC$\/~onions''
in either orientation.  Similar remarks apply to the patterns
showing the octet and dozen with $\sigma^0=0$.

\pagebreak

{\bf 5.~Conjecture B, and other variations and generalizations.}

The problem of determining $\delta(n)$ can be varied and
generalized in many natural ways.  Perhaps the most obvious
is to change the definition of adjacency in~$L$ so $x,y$
are adjacent only if they differ by a unit vector, so each $x\in L$
has 4 neighbors.  Equivalently, if $L$ is identified with $\Z[i]$
then $x,x'$ are adjacent $\iff |x-x'|=1$ (in the usual norm
on~{\bf C}) $\iff x' = x\pm1$ or $x\pm i$.  We then define
$\delta_4(n)$ to be the maximum density of a set $S\subseteq L$
each of whose elements has at most $n$ neighbors in~$L$ relative
to this smaller 4-point neighborhood.  Alternatively we may
vary the lattice; for instance we may replace $L$ by the triangular
lattice of Eisenstein integers
$E=\Z[\rho]$ (with $\rho=e^{2\pi i/3}=(-1+\sqrt{-3})/2$
a cube root of unity), in which each point has six neighbors,
and define $\delta_6(n)$ accordingly.  As it happens all the
$\delta_4$ and $\delta_6$ values are readily obtained using the
methods of~\S2 of this paper:

{\bf Proposition 3.} {\sl The values of $\delta_4$ and $\delta_6$
are given by the following tables:
$$
\begin{array}{c|ccccr}
	n   &  0  &  1  &  2  &  3  & \geqs 4
\\ \hline
\delta_4(n) & 1/2 & 1/2 & 2/3 & 4/5 & 1
\end{array}
$$
$$
\begin{array}{c|ccccccr}
	n   &  0  &  1  &  2  &  3  &  4  &  5  & \geqs6
\\ \hline
\delta_6(n) & 1/3 & 2/5 & 1/2 & 2/3 & 3/4 & 6/7 & 1
\end{array}
$$
These maximal densities are attained by the following periodic sets
(see Appendix~C):\nls
$\delta_4(0)=\delta_4(1)=1/2$ by the checkerboard sublattice
$\{ x\in L: x_1 \equiv x_2 \bmod 2 \}$,\nls
$\delta_4(2)=2/3$ by $\{ x\in L: x_1 \not\equiv x_2 \bmod 3 \}$,\nls
$\delta_4(3)=4/5$ by $\{ x\in L : x_1+2x_2 \not\equiv 0 \bmod 5\}$
(the same set that attains $\delta(6)=4/5$), and\nls
$\delta_4(4)=1$ by~$L$;\nls
$\delta_6(0)=1/3$ by the triangular sublattice $\sqrt{-3}E$,\nls
$\delta_6(1)=2/5$ by the translates of the domino $\{0,1\}$
by \Z-linear combinations of $\sqrt{-3}$ and $2-\rho$,\nls
$\delta_6(2)=1/2$ by $\Z[\sqrt{-3}]$,\nls
$\delta_6(3)=2/3$ by $E - \sqrt{-3} E$,\nls
$\delta_6(4)=3/4$ by $E - 2E$,\nls
$\delta_6(5)=6/7$ by $E - (2-\rho) E$, and\nls
$\delta_6(6)=1$ by~$E$.\nls
Moreover, these are the unique periodic subsets attaining those
densities up to the symmetries of the respective lattices, except
for the cases of $\delta_4(1)$, $\delta_6(2)$ and $\delta_6(4)$.
}

{\em Proof sketch}\/: We readily check that these sets attain the
densities claimed, so we need only prove that those densities
are the largest possible.  First we obtain the bounds
$\delta_4(n)\leqs 4/(8-n)$, $\delta_6(n)\leqs 6/(12-n)$,
attained by a periodic set~$S$\/ if and only if each $x\in S$\/ has
exactly $n$ neighbors in~$S$ and no two elements of $L-S$\/
are adjacent.  This is analogous to Prop.~1 and proved in
exactly the same way; it accounts for all the $\delta_4$'s
excepting $\delta_4(1)$, and for $\delta_6(n)$ for $n\geqs3$.
The same modification of the argument that we used to prove
$\delta(2)=1/2$ also works for $\delta_6(0)=1/3$: if no
two elements of~$S\subset E$\/ are adjacent then each
$y\in E-S$\/ has at most three neighbors in~$S$, so
$S$\/ has density at most $3/(3+6)=1/3$.  The remaining
cases are handled by a Voronoi-cell argument: $\delta_4(1)=1/2$
using the modified $l_\infty$ norm of \S2,~\S3, and $\delta_6(1)=2/5$,
$\delta_6(2)=1/2$ using the Euclidean norm.  (Alternatively we could
have used the Voronoi approach to obtain $\delta_4(n)$ for all~$n$,
and also the inequality $\delta_6(n)\leqs 2/(6-n)$, corresponding to
Prop.~2, which is sharp for $n\leqs3$.)  The uniqueness of the
optimal patterns except for the three cases indicated is checked
as in~\S2; for those three cases, alternative optimal~$S$\/ patterns
are shown in Appendix~C (the first two of which are analogous to
the ``chicken wire'' and ``tv~static'' of Appendix~A).~~$\Qed$

[In particular we see that $\delta_6$ is strictly increasing, while
$\delta_4$, like $\delta$, is not because $\delta_4(0)=\delta_4(1)$.]

A natural further generalization is to lattices in higher dimensions.
Here we have a great choice of neighborhoods and problems --- already
for $\Z^3$ with the $l_\infty$ neighborhood each point of $\Z^3$ is
adjacent to 26 others, so there are 27 maximal densities to obtain,
all but the first and last few of which seem to pose very difficult
problems.  By analogy with the two-dimensional case we expect to
obtain more tractable problems by restricting adjacency to points
differing by a unit vector.  Generalizing to $\Z^k$ for
arbitrary~$k$,\footnote{
  Other generalizations can be considered as well; for isntance
  we could, expanding on a suggestion of Greg Kuperberg, use any
  root lattice, with $x,x'$ adjacent if $x-x'$ is a root vector,
  or perhaps a short root vector.}
we thus define $\Delta_k(n)$ to be the maximal upper density of a subset
$S\in\Z^k$ such that for each $x\in S$\/ there are at most $n$ points
$x'\in S$\/ at (Euclidean) distance~$1$ from~$x$.  (Our usual
cut-and-paste argument shows that ``upper'' is superfluous and
the maximum density is attained by some~$S$.)  Note that
$\Delta_2$ is what we called $\delta_4$ in Prop.~3.  Preliminary
work suggests that these $\Delta_k(n)$ do in fact constitute
a promising generalization.  We first collect some easy results:

{\bf Proposition 4.}  {\sl For all $k,n$ we have
$\Delta_{k+1}(n) \leqs \Delta_k(n) \leqs \Delta_k(n+1)$,
i.e.\ $\Delta_k(n)$ is nonincreasing as a function of~$k$
and nondecreasing as a function of~$n$.  Also
$\Delta_{mk}(mn)\geqs\Delta_k(n)$ for any positive integer~$m$,
and $\Delta_k(n)\leqs 2k/(4k-n)$, with equality at least when
$n=2k$ or $(2k-n)|k$ or $(2k+n)/(2k-n)$ is a power of~$2$.}

[This last condition may be stated equivalently as follows:
either $n=2k$ or $1-(n/2k)$ is $1/2l$ or $1/(2^l-1)$ for
some positive integer~$l$.  In particular $n=0$ is allowed.]

{\em Proof}\/:  That $\Delta_k(n)\leqs \Delta_k(n+1)$ is clear.
If $S\subset\Z^{k+1}$ has maximal degree $\leqs n$ then so do
its $k$-dimensional slices, and if moreover $S$\/ has density $\Delta$
then it has slices of density $>\Delta-\eps$ for each $\eps>0$;
thus $\Delta_{k+1}(n) \leqs \Delta_k(n)$.  Moreover the subset
$$
S^{[m]} := \{ (x^{(1)}, x^{(2)}, \ldots, x^{(m)}) :
  x^{(i)} \in \Z^k, \sum_{i=1}^m x^{(i)} \in S \}
$$
of $\Z^{mk}$ has maximal degree $\leqs mn$ and density\footnote{
  This is actually tricky to prove with our definition of the
  density.  However if $\Delta=\Delta_k(n)$ then translates
  of {\sf B}$_r$ in $S$\/ have density $\Delta\pm O(1/r)$, from which
  the fact that $S^{[m]}$ has density~$\Delta$ quickly follows.
  }
$\Delta$, whence $\Delta_{mk}(mn)\geqs\Delta_k(n)$.  Finally
$\Delta_k(n)\leqs 2k/(4k-n)$ is proved as in Prop.~1, and since
this upper bound depends only on $n/2k$ equality holds for
$\Delta_{mk}(mn)$ if it holds for $\Delta_k(n)$.  Thus to
complete the proof the this Proposition we need only show that
the bounds on $\Delta_1(1)$, $\Delta_k(2k-1)$, and
$\Delta_{2^l-1}(2^{l+1}-4)$ are attained.  The first is
clear (let $S=\Z$), and the other two are attained by sublattice
complements:
$$
\{ x \in \R^k: \sum_{i} i x_i \not\equiv 0 \bmod 2k+1 \}
$$
has degree $2k-1$ and density $2k/(2k+1)=2k/(4k-n)$,
and the complement of the index-$2^l$ lattice
consisting of $x\in\Z^{2^l-1}$ such that $x\bmod 2$ is in
the (perfect one-error-correcting) extended Hamming code
with parameters $[2^l-1,2^l-l-1,3]$ (see \cite[p.23--25]{TECC})
has degree $2^{l+1}-4$ and density $(2^l-1)/2^l$.~~\Qed

In particular $\Delta_k(0)=1/2$ for all $k$, attained by the
checkerboard lattice
$ 2 D^*_k = \{ x \in \Z^k: \sum_{i=1}^k x_i \equiv 0 \bmod 2 \}$;
and $\Delta_3(4)=3/4$, attained by the complement of the bcc lattice
$\{ x \in \Z^3: x_1 \equiv x_2 \equiv x_3 \bmod 2 \}$.

{\bf Corollary} {\sl We have $\Delta_k(1)=1/2$ for all $k\geqs2$;
thus $\Delta_k$ is not a strictly increasing function of~$n$
for any $k\geqs2$.}

{\em Proof}\/: Prop.~3 gives $\Delta_2(1)=\delta_4(1)=1/2$.
By Prop.~4, for any $k\geqs 2$
$$
\frac12 = \Delta_k(0) \leqs \Delta_k(1) \leqs \Delta_2(1) = \frac12,
$$
so equality must hold throughout.~~\Qed

Of course $\Delta_1(1)=2/3$, which is also a special case of Prop.~4.
We now have all values of $\Delta_k(n)$ for $k\leqs3$ except for
$\Delta_3(2)$.  By Prop.~4 this is bounded above by $3/5$, but
equality cannot hold: if for some $S\in\Z^3$ every $x\in S$\/
has exactly two neighbors in~$S$\/ then every $y\in \Z^3 - S$
has at least two neighbors in~$\Z^3 - S$, so $S$\/ has density
at most $1/2$.  In fact I have not been able to find a single
$S\in\Z^3$ of maximal degree~2 and density $>1/2$, or for that
matter any $S\in\Z^k$ of maximal degree $<k$ and density
$>1/2=\Delta_k(0)$.  I thus propose, perhaps rashly, that

{\bf Conjecture B.}  {\sl If $0 \leqs n < k$ then $\Delta_k(n)=1/2$.}

If true, this has the surprising consequence that, for each $k$,
the a priori nondecreasing function $\Delta_k(n)$ is actually
constant for $0\leqs n<k$, then jumps from $1/2$ to $2/3$ at $k=n$.
This would still leave open the question of $\Delta_k(n)$ for
$n<k<2n$; by Prop.~4 this is $2k/(4k-n)$ for $k\leqs 4$ with the
possible exception of $\Delta_4(5)$, but it seems unlikely that
$\Delta_4(5)=8/11$.

To generalize this problem further we could expand the neighborhood
further and allow different weights, with the weight of a ``neighbor''
$x'$ of~$x$ depending only on $x'-x$, not necessarily respecting
any lattice symmetries.  That is, we fix a positive integer~$k$,
a real number~$c$, and a function $w: \Lambda \ra \R$ supported
on a finite subset of a lattice $\Lambda\cong\Z^k$, and then
ask for the maximum density $\DD_w(c)$ of a subset~$S$\/
of~$\Z^k$ such that
$$
x \in S \Rightarrow \sum_{x'\in S} w(x'-x) \leqs c.
\eqno{(*)}
$$
(For example $\Delta_k(n)=\DD_w(n)$ where $w(v)$ is~$1$
if $|v|=1$ and $0$ otherwise.  Note that in general the $w(v)$
are not required to be positive.)  At this level of generality
it is probably not reasonable to expect to compute every $\DD_w(c)$.
But a more basic problem presents itself.  Since all our
conditions on~$S$\/ are local, the same cut-and-paste
argument we used for $\delta(n)$ etc.\ shows that
the maximum density is attained, and equals the maximum
upper density.  Moreover if $\Delta$ is that maximum
then (*) is satisfied by periodic sets of density $>\Delta-\eps$
for each $\eps>0$.  However, in each of the cases where
we have been able to determine a maximum density it was
actually attained (non-uniquely) by a periodic set, and
it is not a priori clear whether this should be true in
general.  We thus ask:

{\bf Question.}  {\sl Is $\DD_w(c)$ [or even $\Delta_k(n)$
or $\delta(n)$] necessarily attained by a periodic set?
In particular is every such maximum density rational?}

I do not call this a Conjecture because I do not have
particularly good reason to either believe it is true
or suspect it is false, though
either a proof or a disproof would be of interest.
One might propose an even more far-reaching generalization,
which accommodates local conditions not of the form~(*):
choose some $r>0$ and a family $\FF$ of subsets of $\B_r\subset\Lambda$,
and let $\DD(\FF)$ be the maximum density of $S\subseteq\Lambda$
satisfying
$$
x \in S \Rightarrow \{x'-x: x'\in S\} \cap \B_r \in \FF.
$$
But this is probably too ambitious a generalization,
because it should be possible to find $\Delta,\FF$
such that the question ``does $\DD(\FF)\geqs\Delta$\,?''\
encodes the general tiling problem, and it is known
(see e.g.~\cite{Robinson})
that the solution of the general tiling problem is not
computable.  It should also be possible to find in the same
way $\FF$ such that $\DD(\FF)$ is an arbitrary computable
irrational, and/or attained only by subsets of~$\Lambda$
that are not even almost periodic.

From these clouds of abstraction we return home to our
original motivation in the game of Life for a final
generalization of a rather different flavor.  Recall
that a still-Life is a fixed point of the map
$\lambda: \SS\ra\SS$.  More generally, an {\em oscillator}
is a periodic point\footnote{
  Warning: this is ``time periodicity'' as opposed to the
  ``space periodicity'' seen thus far.
  }
of~$\lambda$, i.e.\ an $S\subset L$
such that $\lambda^p S = S$ for some $p>0$; the least
such~$p$ is the {\em period} of~$S$.  For instance,
a still-Life is an oscillator of period~$1$; the ``blinker''
$\{ (x_1,0) : |x_1| \leqs 1 \}$ is the smallest
oscillator of period~$2$ (see Figure~6a; this is also
the smallest nonempty oscillator of any period).
Two examples of infinite oscillators are the period-2
``Venetian blinds'' $\{ x: x_1 \equiv 0 {\rm\ or\ }1 \bmod 4 \}$,
and a period-6 analogue $\{ x: x_1 \equiv 0 {\rm\ or\ }1 \bmod 8 \}$.
There are oscillators of period~$p$ for any $p>0$;
one of many constructions proving this is the ``lightspeed wire''
shown in Fig.~6b, which easily generalizes from
the shown $p=7$ to arbitrary~$p\geqs 5$.\footnote{
  It is likely that {\em finite} oscillators of any period exist,
  though this is still an open problem despite considerable
  recent progress.  It has been known for some time (unpublished
  work of Conway and his students) that there exists
  an integer $p_0$ such that the ``universal constructors'' described
  in~\cite{WW} yield a finite oscillator of period~$p$ for each
  $p\geqs p_0$, but such an oscillator has never been exhibited,
  and even a specific value of $p_0$ was never established.
  Only this year it was shown~\cite{BC}, using a completely different
  approach, how oscillators of all periods $\geqs 58$ may be
  constructed explicitly.  As of this writing the only $p$ still
  unknown are 19, 23, 27, 31, 37, 38, 41, 43, 49, 53, and~57.
     [But see Note Added in Proof next page!]
  Also the only oscillators known of periods $p=33$, 34, 39, 51
  are trivially interacting pairs of oscillators whose periods
  are the prime factors of~$p$.  For all other~$p$ a finite oscillator
  of genuine period~$p$ is known.  In particular oscillators of period~3
  and~4 have been known for some time, so indeed every $p$ is realized
  by some possibly infinite oscillator.
  }
Note that two of the phases $\{ x: x_1 \not\equiv 4,5 \bmod 8\}$,
$\{ x: x_1 \not\equiv 0,1 \bmod 8\}$ of our period-6 oscillator
have density $3/4$; thus we cannot hope to prove that any oscillator
has density $\leqs1/2$.  However for an oscillator a more natural
invariant is the average density of its phases, and that average
density is never observed to exceed $1/2$, though $1/2$ is attained
by some oscillators of period $>1$, such as ``Venetian blinds''$\!$.
We thus ask whether in fact every Life oscillator has average
density $\leqs1/2$, and also (knowing no counterexample) whether
each phase of a period-2 oscillator has density $\leqs1/2$.
While we have seen that an oscillator may have density $>1/2$ in
some of its phases, we do not know what is the largest density 
possible in a single oscillator: can it come arbitrarily close
to~$1$?  Can it even exceed~$3/4$?  It will be noted that these
problems are either special cases or closely related to special
cases of the $\DD(\FF)$ problem in~2 or~3 dimensions, but perhaps
these special cases are more tractable than the problems with
arbitrary~$\FF$.
\def\b{\makebox[.3cm]{$\bullet$}}
\def\0{\makebox[.3cm]{$\cdot$}}
\def\:{\makebox[.3cm]{}}
\renewcommand{\arraystretch}{.8}
$$
\sq{
\0\0\0\0\0\makebox[1.2cm]{}\0\0\0\0\0\\
\0\0\b\0\0\makebox[1.2cm]{}\0\0\0\0\0\\
\0\0\b\0\0\makebox[1.2cm]{$\longleftrightarrow$}\0\b\b\b\0\\
\0\0\b\0\0\makebox[1.2cm]{}\0\0\0\0\0\\
\0\0\0\0\0\makebox[1.2cm]{}\0\0\0\0\0\\
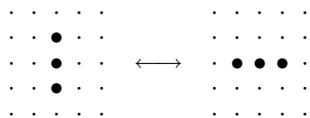
\\ \mbox{\large Fig.~6a: blinker (period 2)}
}
\qquad
\sq{
\0\b\b\0\0\b\b\0\0\b\b\0\0\b\b\0\0\b\b\0\\
\0\b\b\0\0\b\b\0\0\b\b\0\0\b\b\0\0\b\b\0\\
\0\0\0\0\0\0\0\0\0\0\0\0\0\0\0\0\0\0\0\0\\
\b\b\b\b\b\b\b\b\b\b\b\b\b\b\b\b\b\b\b\b\\
\0\0\b\0\0\0\0\0\b\0\0\0\0\0\0\0\b\0\0\0\\
\0\0\0\b\b\b\0\0\0\0\b\b\b\0\0\0\0\b\b\b\\
\0\b\0\0\0\0\0\0\0\b\0\0\0\0\0\b\0\0\0\0\\
\b\b\b\b\b\b\b\b\b\b\b\b\b\b\b\b\b\b\b\b\\
\0\0\0\0\0\0\0\0\0\0\0\0\0\0\0\0\0\0\0\0\\
\0\b\b\0\0\b\b\0\0\b\b\0\0\b\b\0\0\b\b\0\\
\0\b\b\0\0\b\b\0\0\b\b\0\0\b\b\0\0\b\b\0\\
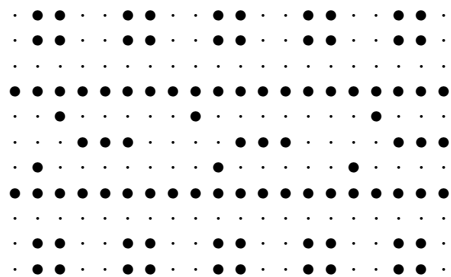
\\ \mbox{\large Fig.~6b: lightspeed wire (period 7)}
}
$$

\pagebreak

{\bf Acknowledgements.}  Thanks to Dan Hoey and other members
of the Life mailing list for correcting an earlier version
of my proof of Conjecture~A and, together with Greg Kuperberg and
Marcus Moore, for helpful e-mail discussions; thanks also for
correcting an earlier draft of this manuscript to Dan Hoey (again)
and to Dean Hickerson.

This work was made possible in part by funding from the
National Science Foundation and the Packard Foundation.

{\bf Added in Proof:} On 18 December 1997, Dietrich Leithner constructed
a finite Life pattern of period 57 by adapting the methods of~\cite{BC}.
Thus 57 can be struck from the list of unknown periods in Footnote 5;
finite oscillators of all periods $\geqs 54$ are now known.

{\bf Appendix A: A gallery of density-$1/2$ still Lifes.}

\vspace*{3ex}

We exhibit and name several periodic or almost-periodic still-Life
patterns attaining the maximal density 1/2 in various ways ---
indeed in all the ways known to us.  The names ``chicken wire''
and ``onion rings'' have been in use in the Life community for
some time, and I have not been able to determine their source.
``Onion bulbs'' is a variant on ``onion rings''$\!$, while
``tv static'' and ``square waves'' are descriptive neologisms.
The Hoey and Moore patterns are named after their discoverers;
concerning ``DCC onions''$\!$, ``octets'' and ``dozens''
see~\S4.

\vspace*{5ex}

\centerline{
$
\sq{
\0\0\b\b\0\0\b\b\b\0\0\0\0\0\b\b\0\0\0\\
\b\b\0\0\b\b\0\0\0\b\b\b\b\b\0\0\b\b\b\\
\0\0\b\b\0\0\b\b\b\0\0\0\0\0\b\b\0\0\0\\
\b\b\0\0\b\b\0\0\0\b\b\b\b\b\0\0\b\b\b\\
\0\0\b\b\0\0\b\b\b\0\0\0\0\0\b\b\0\0\0\\
\b\b\0\0\b\b\0\0\0\b\b\b\b\b\0\0\b\b\b\\
\0\0\b\b\0\0\b\b\b\0\0\0\0\0\b\b\0\0\0\\
\b\b\0\0\b\b\0\0\0\b\b\b\b\b\0\0\b\b\b\\
\0\0\b\b\0\0\b\b\b\0\0\0\0\0\b\b\0\0\0\\
\b\b\0\0\b\b\0\0\0\b\b\b\b\b\0\0\b\b\b\\
\0\0\b\b\0\0\b\b\b\0\0\0\0\0\b\b\0\0\0\\
\b\b\0\0\b\b\0\0\0\b\b\b\b\b\0\0\b\b\b\\
\0\0\b\b\0\0\b\b\b\0\0\0\0\0\b\b\0\0\0\\
\b\b\0\0\b\b\0\0\0\b\b\b\b\b\0\0\b\b\b\\
\0\0\b\b\0\0\b\b\b\0\0\0\0\0\b\b\0\0\0\\
\b\b\0\0\b\b\0\0\0\b\b\b\b\b\0\0\b\b\b\\
\0\0\b\b\0\0\b\b\b\0\0\0\0\0\b\b\0\0\0\\
\b\b\0\0\b\b\0\0\0\b\b\b\b\b\0\0\b\b\b\\
\0\0\b\b\0\0\b\b\b\0\0\0\0\0\b\b\0\0\0\\
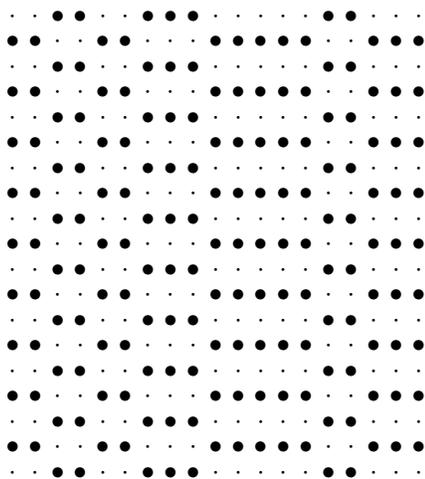
\\ \mbox{\large\rm chicken wire (still LiFe)}
}
\qquad
\sq{
\0\0\0\0\b\b\b\b\0\0\0\0\0\b\b\b\b\0\0\\
\b\b\b\0\0\0\0\b\b\b\b\b\0\0\0\0\b\b\b\\
\0\0\b\b\b\b\0\0\0\0\0\b\b\b\b\0\0\0\0\\
\b\0\0\0\0\b\b\b\b\b\0\0\0\0\b\b\b\b\0\\
\b\b\b\b\0\0\0\0\0\b\b\b\b\0\0\0\0\b\0\\
\0\0\0\b\b\b\b\b\0\0\0\0\b\b\b\b\0\b\0\\
\b\b\0\0\0\0\0\b\b\b\b\0\0\0\0\b\0\b\0\\
\0\b\b\b\b\b\0\0\0\0\b\b\b\b\0\b\0\b\b\\
\0\0\0\0\0\b\b\b\b\0\0\0\0\b\0\b\0\0\b\\
\b\b\b\b\0\0\0\0\b\b\b\b\0\b\0\b\b\0\b\\
\0\0\0\b\b\b\b\0\0\0\0\b\0\b\0\0\b\0\b\\
\b\b\0\0\0\0\b\b\b\b\0\b\0\b\b\0\b\0\0\\
\0\b\b\b\b\0\0\0\0\b\0\b\0\0\b\0\b\b\b\\
\0\0\0\0\b\b\b\b\0\b\0\b\b\0\b\0\0\0\0\\
\b\b\b\0\0\0\0\b\0\b\0\0\b\0\b\b\b\b\0\\
\0\0\b\b\b\b\0\b\0\b\b\0\b\0\0\0\0\b\b\\
\b\0\0\0\0\b\0\b\0\0\b\0\b\b\b\b\0\0\0\\
\b\b\b\b\0\b\0\b\b\0\b\0\0\0\0\b\b\b\b\\
\0\0\0\b\0\b\0\0\b\0\b\b\b\b\0\0\0\0\b\\
\\ \mbox{\large\rm tv\ static}
}
$
}
\centerline{
$
\sq{
\b\0\b\0\b\0\0\b\0\b\0\b\b\0\b\0\b\0\0\\
\0\0\b\0\b\b\b\b\0\b\0\0\0\0\b\0\b\b\b\\
\b\b\b\0\0\0\0\0\0\b\b\b\b\b\b\0\0\0\0\\
\0\0\0\b\b\b\b\b\b\0\0\0\0\0\0\b\b\b\b\\
\b\b\0\b\0\0\0\0\b\0\b\b\b\b\0\b\0\0\0\\
\0\b\0\b\0\b\b\0\b\0\b\0\0\b\0\b\0\b\b\\
\0\b\0\b\0\b\b\0\b\0\b\0\0\b\0\b\0\b\b\\
\b\b\0\b\0\0\0\0\b\0\b\b\b\b\0\b\0\0\0\\
\0\0\0\b\b\b\b\b\b\0\0\0\0\0\0\b\b\b\b\\
\b\b\b\0\0\0\0\0\0\b\b\b\b\b\b\0\0\0\0\\
\0\0\b\0\b\b\b\b\0\b\0\0\0\0\b\0\b\b\b\\
\b\0\b\0\b\0\0\b\0\b\0\b\b\0\b\0\b\0\0\\
\b\0\b\0\b\0\0\b\0\b\0\b\b\0\b\0\b\0\0\\
\0\0\b\0\b\b\b\b\0\b\0\0\0\0\b\0\b\b\b\\
\b\b\b\0\0\0\0\0\0\b\b\b\b\b\b\0\0\0\0\\
\0\0\0\b\b\b\b\b\b\0\0\0\0\0\0\b\b\b\b\\
\b\b\0\b\0\0\0\0\b\0\b\b\b\b\0\b\0\0\0\\
\0\b\0\b\0\b\b\0\b\0\b\0\0\b\0\b\0\b\b\\
\0\b\0\b\0\b\b\0\b\0\b\0\0\b\0\b\0\b\b\\
\b\b\0\b\0\0\0\0\b\0\b\b\b\b\0\b\0\0\0\\
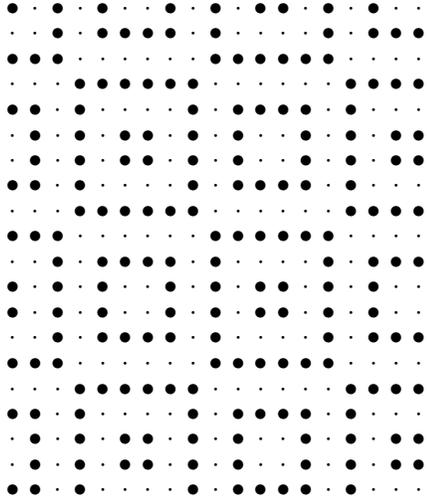
\\ \mbox{\large\rm onion rings (order 3)}
}
\qquad
\sq{
\b\0\b\b\0\b\0\0\b\0\b\b\0\b\0\0\b\0\b\\
\b\0\0\0\0\b\b\b\b\0\0\0\0\b\b\b\b\0\0\\
\0\b\b\b\b\0\0\0\0\b\b\b\b\0\0\0\0\b\b\\
\b\0\0\0\0\b\b\b\b\0\0\0\0\b\b\b\b\0\0\\
\0\b\b\b\b\0\0\0\0\b\b\b\b\0\0\0\0\b\b\\
\0\b\0\0\b\0\b\b\0\b\0\0\b\0\b\b\0\b\0\\
\0\b\0\0\b\0\b\b\0\b\0\0\b\0\b\b\0\b\0\\
\0\b\b\b\b\0\0\0\0\b\b\b\b\0\0\0\0\b\b\\
\b\0\0\0\0\b\b\b\b\0\0\0\0\b\b\b\b\0\0\\
\0\b\b\b\b\0\0\0\0\b\b\b\b\0\0\0\0\b\b\\
\b\0\0\0\0\b\b\b\b\0\0\0\0\b\b\b\b\0\0\\
\0\b\b\b\b\0\0\0\0\b\b\b\b\0\0\0\0\b\b\\
\0\b\0\0\b\0\b\b\0\b\0\0\b\0\b\b\0\b\0\\
\0\b\0\0\b\0\b\b\0\b\0\0\b\0\b\b\0\b\0\\
\0\b\b\b\b\0\0\0\0\b\b\b\b\0\0\0\0\b\b\\
\b\0\0\0\0\b\b\b\b\0\0\0\0\b\b\b\b\0\0\\
\0\b\b\b\b\0\0\0\0\b\b\b\b\0\0\0\0\b\b\\
\b\0\0\0\0\b\b\b\b\0\0\0\0\b\b\b\b\0\0\\
\b\0\b\b\0\b\0\0\b\0\b\b\0\b\0\0\b\0\b\\
\b\0\b\b\0\b\0\0\b\0\b\b\0\b\0\0\b\0\b\\
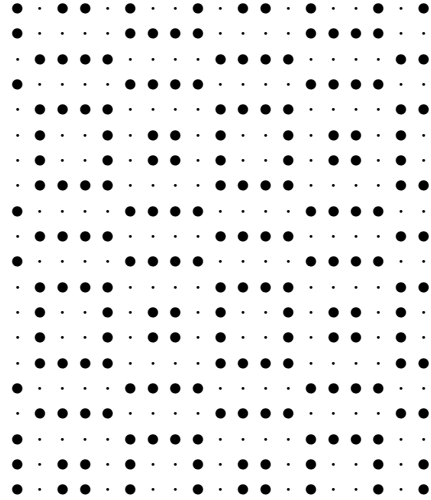
\\ \mbox{\large\rm onion bulbs (order 2)}
}
$
}

\vspace*{5ex}

\centerline{
$
\sq{
\0\b\0\0\0\b\0\b\0\0\0\b\0\b\0\0\0\b\0\b\\
\0\0\b\b\b\b\0\0\b\b\b\b\0\0\b\b\b\b\0\0\\
\b\0\b\0\0\0\b\0\b\0\0\0\b\0\b\0\0\0\b\0\\
\b\0\0\b\b\b\b\0\0\b\b\b\b\0\0\b\b\b\b\0\\
\0\b\0\b\0\0\0\b\0\b\0\0\0\b\0\b\0\0\0\b\\
\b\b\0\0\b\b\b\b\0\0\b\b\b\b\0\0\b\b\b\b\\
\0\0\b\0\b\0\0\0\b\0\b\0\0\0\b\0\b\0\0\0\\
\b\b\b\0\0\b\b\b\b\0\0\b\b\b\b\0\0\b\b\b\\
\0\0\0\b\0\b\0\0\0\b\0\b\0\0\0\b\0\b\0\0\\
\b\b\b\b\0\0\b\b\b\b\0\0\b\b\b\b\0\0\b\b\\
\0\0\0\0\b\0\b\0\0\0\b\0\b\0\0\0\b\0\b\0\\
\b\b\b\b\b\0\0\b\b\b\b\0\0\b\b\b\b\0\0\b\\
\b\0\0\0\0\b\0\b\0\0\0\b\0\b\0\0\0\b\0\b\\
\0\b\b\b\b\b\0\0\b\b\b\b\0\0\b\b\b\b\0\0\\
\0\b\0\0\0\0\b\0\b\0\0\0\b\0\b\0\0\0\b\0\\
\0\0\b\b\b\b\b\0\0\b\b\b\b\0\0\b\b\b\b\0\\
\b\0\b\0\0\0\0\b\0\b\0\0\0\b\0\b\0\0\0\b\\
\b\0\0\b\b\b\b\b\0\0\b\b\b\b\0\0\b\b\b\b\\
\0\b\0\b\0\0\0\0\b\0\b\0\0\0\b\0\b\0\0\0\\
\b\b\0\0\b\b\b\b\b\0\0\b\b\b\b\0\0\b\b\b\\
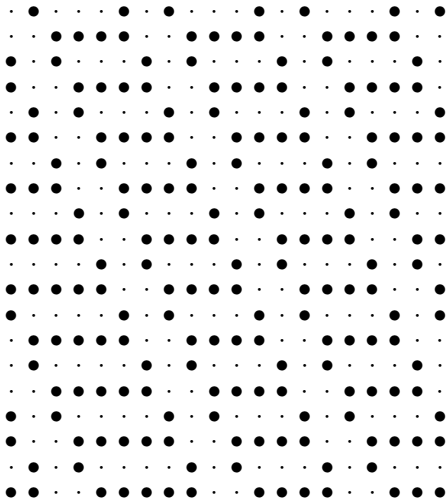
\\ \mbox{\large Hoey I (orders 1, 2)}
}
\qquad
\sq{
\0\0\0\0\b\0\0\0\0\b\0\0\0\0\b\0\0\0\0\\
\b\b\b\b\0\b\b\b\b\0\b\b\b\b\0\b\b\b\b\\
\b\0\0\b\0\b\0\0\b\0\b\0\0\b\0\b\0\0\b\\
\b\0\b\0\0\b\0\b\0\0\b\0\b\0\0\b\0\b\0\\
\b\0\b\b\b\b\0\b\b\b\b\0\b\b\b\b\0\b\b\\
\0\b\0\0\0\0\b\0\0\0\0\b\0\0\0\0\b\0\0\\
\0\b\0\b\b\0\b\0\b\b\0\b\0\b\b\0\b\0\b\\
\0\b\0\b\b\0\b\0\b\b\0\b\0\b\b\0\b\0\b\\
\0\b\0\0\0\0\b\0\0\0\0\b\0\0\0\0\b\0\0\\
\b\0\b\b\b\b\0\b\b\b\b\0\b\b\b\b\0\b\b\\
\b\0\b\0\0\b\0\b\0\0\b\0\b\0\0\b\0\b\0\\
\0\0\b\0\b\0\0\b\0\b\0\0\b\0\b\0\0\b\0\\
\b\b\b\0\b\b\b\b\0\b\b\b\b\0\b\b\b\b\0\\
\0\0\0\b\0\0\0\0\b\0\0\0\0\b\0\0\0\0\b\\
\b\b\0\b\0\b\b\0\b\0\b\b\0\b\0\b\b\0\b\\
\b\b\0\b\0\b\b\0\b\0\b\b\0\b\0\b\b\0\b\\
\0\0\0\b\0\0\0\0\b\0\0\0\0\b\0\0\0\0\b\\
\b\b\b\0\b\b\b\b\0\b\b\b\b\0\b\b\b\b\0\\
\0\0\b\0\b\0\0\b\0\b\0\0\b\0\b\0\0\b\0\\
\0\b\0\0\b\0\b\0\0\b\0\b\0\0\b\0\b\0\0\\
\\ \mbox{\large\rm Hoey II}
}
$
}

\pagebreak

\vspace*{-5ex}

\centerline{
$
\sq{
\b\b\b\b\b\b\b\b\b\b\b\b\b\b\b\b\b\b\b\\
\b\0\0\0\0\0\0\0\b\0\0\0\0\0\0\0\b\0\0\\
\0\0\b\b\b\b\b\0\0\0\b\b\b\b\b\0\0\0\b\\
\b\b\b\0\0\0\b\b\b\b\b\0\0\0\b\b\b\b\b\\
\0\0\0\0\b\0\0\0\0\0\0\0\b\0\0\0\0\0\0\\
\b\b\b\b\b\b\b\b\b\b\b\b\b\b\b\b\b\b\b\\
\0\0\0\0\0\0\0\0\0\0\0\0\0\0\0\0\0\0\0\\
\b\b\b\b\b\b\b\b\b\b\b\b\b\b\b\b\b\b\b\\
\b\0\0\0\b\0\0\0\b\0\0\0\b\0\0\0\b\0\0\\
\0\0\b\0\0\0\b\0\0\0\b\0\0\0\b\0\0\0\b\\
\b\b\b\b\b\b\b\b\b\b\b\b\b\b\b\b\b\b\b\\
\0\0\0\0\0\0\0\0\0\0\0\0\0\0\0\0\0\0\0\\
\b\b\b\b\b\b\b\b\b\b\b\b\b\b\b\b\b\b\b\\
\b\0\0\0\0\0\0\0\b\0\0\0\0\0\0\0\b\0\0\\
\0\0\b\b\b\b\b\0\0\0\b\b\b\b\b\0\0\0\b\\
\b\b\b\0\0\0\b\b\b\b\b\0\0\0\b\b\b\b\b\\
\0\0\0\0\b\0\0\0\0\0\0\0\b\0\0\0\0\0\0\\
\b\b\b\b\b\b\b\b\b\b\b\b\b\b\b\b\b\b\b\\
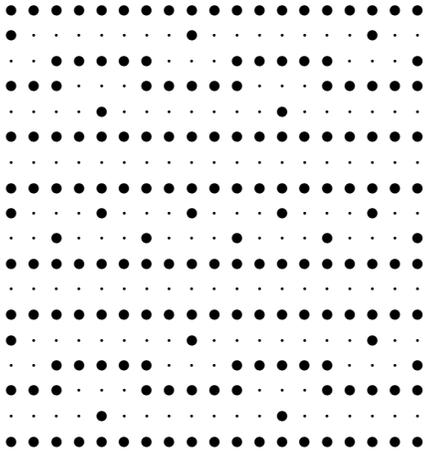
\\ \mbox{\large\rm square waves (orders 1, 2)}
}
\qquad
\sq{
\b\b\b\0\0\b\0\b\0\0\b\b\b\0\b\0\b\0\0\0\b\0\b\\
\b\0\0\b\b\b\0\b\0\b\0\0\0\b\0\b\0\b\b\b\0\0\b\\
\b\0\b\0\0\0\b\0\b\0\b\b\b\0\0\b\0\b\0\0\b\b\b\\
\0\b\0\b\b\b\0\0\b\0\b\0\0\b\b\b\0\b\0\b\0\0\0\\
\0\b\0\b\0\0\b\b\b\0\b\0\b\0\0\0\b\0\b\0\b\b\b\\
\b\b\0\b\0\b\0\0\0\b\0\b\0\b\b\b\0\0\b\0\b\0\0\\
\0\0\b\0\b\0\b\b\b\0\0\b\0\b\0\0\b\b\b\0\b\0\b\\
\b\b\0\0\b\0\b\0\0\b\b\b\0\b\0\b\0\0\0\b\0\b\0\\
\0\0\b\b\b\0\b\0\b\0\0\0\b\0\b\0\b\b\b\0\0\b\0\\
\0\b\0\0\0\b\0\b\0\b\b\b\0\0\b\0\b\0\0\b\b\b\0\\
\b\0\b\b\b\0\0\b\0\b\0\0\b\b\b\0\b\0\b\0\0\0\b\\
\b\0\b\0\0\b\b\b\0\b\0\b\0\0\0\b\0\b\0\b\b\b\0\\
\b\0\b\0\b\0\0\0\b\0\b\0\b\b\b\0\0\b\0\b\0\0\b\\
\0\b\0\b\0\b\b\b\0\0\b\0\b\0\0\b\b\b\0\b\0\b\0\\
\b\0\0\b\0\b\0\0\b\b\b\0\b\0\b\0\0\0\b\0\b\0\b\\
\0\b\b\b\0\b\0\b\0\0\0\b\0\b\0\b\b\b\0\0\b\0\b\\
\b\0\0\0\b\0\b\0\b\b\b\0\0\b\0\b\0\0\b\b\b\0\b\\
\0\b\b\b\0\0\b\0\b\0\0\b\b\b\0\b\0\b\0\0\0\b\0\\
\\ \mbox{\large\rm Moore}
}
$
}

\vspace*{8ex}

\centerline{
$
\sq{
\0\0\0\b\0\0\0\0\b\0\0\0\0\b\0\0\0\0\b\0\0\0 \\
\b\b\b\0\b\b\b\b\0\b\b\b\b\0\b\b\b\b\0\b\b\b \\
\0\0\b\0\b\0\0\b\0\b\0\0\b\0\b\0\0\b\0\b\0\0 \\
\0\0\b\0\b\0\0\b\0\b\0\0\b\0\b\0\0\b\0\b\0\0 \\
\b\b\b\0\b\b\b\b\0\b\b\b\b\0\b\b\b\b\0\b\b\b \\
\0\0\0\b\0\0\0\0\b\0\0\0\0\b\0\0\0\0\b\0\0\0 \\
\b\b\b\0\0\b\b\b\0\0\b\b\b\0\0\b\b\b\0\0\b\b \\
\b\0\0\b\b\b\0\0\b\b\b\0\0\b\b\b\0\0\b\b\b\0 \\
\0\0\b\0\0\0\0\b\0\0\0\0\b\0\0\0\0\b\0\0\0\0 \\
\b\b\0\b\b\b\b\0\b\b\b\b\0\b\b\b\b\0\b\b\b\b \\
\0\b\0\b\0\0\b\0\b\0\0\b\0\b\0\0\b\0\b\0\0\b \\
\0\b\0\b\0\0\b\0\b\0\0\b\0\b\0\0\b\0\b\0\0\b \\
\b\b\0\b\b\b\b\0\b\b\b\b\0\b\b\b\b\0\b\b\b\b \\
\0\0\b\0\0\0\0\b\0\0\0\0\b\0\0\0\0\b\0\0\0\0 \\
\b\b\0\0\b\b\b\0\0\b\b\b\0\0\b\b\b\0\0\b\b\b \\
\0\0\b\b\b\0\0\b\b\b\0\0\b\b\b\0\0\b\b\b\0\0 \\
\b\0\0\0\0\b\0\0\0\0\b\0\0\0\0\b\0\0\0\0\b\0 \\
\0\b\b\b\b\0\b\b\b\b\0\b\b\b\b\0\b\b\b\b\0\b \\
\0\b\0\0\b\0\b\0\0\b\0\b\0\0\b\0\b\0\0\b\0\b \\
\\ \mbox{\large $DCC$ onions}
}
\qquad
\sq{
\b\0\b\b\b\0\b\0\b\0\b\b\b\0\b\0\b\0\b\b\b\0\b \\
\b\0\0\0\0\0\b\0\b\0\0\0\0\0\b\0\b\0\0\0\0\0\b \\
\b\b\b\b\b\b\b\0\b\b\b\b\b\b\b\0\b\b\b\b\b\b\b \\
\0\0\0\0\0\0\0\b\0\0\0\0\0\0\0\b\0\0\0\0\0\0\0 \\
\0\b\b\b\b\b\0\b\0\b\b\b\b\b\0\b\0\b\b\b\b\b\0 \\
\0\b\0\0\0\b\0\b\0\b\0\0\0\b\0\b\0\b\0\0\0\b\0 \\
\0\b\0\b\0\b\0\b\0\b\0\b\0\b\0\b\0\b\0\b\0\b\0 \\
\b\0\b\b\b\0\b\0\b\0\b\b\b\0\b\0\b\0\b\b\b\0\b \\
\b\0\0\0\0\0\b\0\b\0\0\0\0\0\b\0\b\0\0\0\0\0\b \\
\b\b\b\b\b\b\b\0\b\b\b\b\b\b\b\0\b\b\b\b\b\b\b \\
\0\0\0\0\0\0\0\b\0\0\0\0\0\0\0\b\0\0\0\0\0\0\0 \\
\0\b\b\b\b\b\0\b\0\b\b\b\b\b\0\b\0\b\b\b\b\b\0 \\
\0\b\0\0\0\b\0\b\0\b\0\0\0\b\0\b\0\b\0\0\0\b\0 \\
\0\b\0\b\0\b\0\b\0\b\0\b\0\b\0\b\0\b\0\b\0\b\0 \\
\b\0\b\b\b\0\b\0\b\0\b\b\b\0\b\0\b\0\b\b\b\0\b \\
\b\0\0\0\0\0\b\0\b\0\0\0\0\0\b\0\b\0\0\0\0\0\b \\
\b\b\b\b\b\b\b\0\b\b\b\b\b\b\b\0\b\b\b\b\b\b\b \\
\0\0\0\0\0\0\0\b\0\0\0\0\0\0\0\b\0\0\0\0\0\0\0 \\
\0\b\b\b\b\b\0\b\0\b\b\b\b\b\0\b\0\b\b\b\b\b\0 \\
\0\b\0\0\0\b\0\b\0\b\0\0\0\b\0\b\0\b\0\0\0\b\0 \\
\\ \mbox{\large $AB^5F'F''$ octets}
}
$
}

\centerline{
$
\sq{
\b\0\0\0\0\0\0\0\b\0\0\0\0\0\0\0\b\0\0\0\0\0\0\0\b\0\0\0\0\0\0\0\b \\
\0\b\b\b\b\b\b\b\b\b\b\b\b\b\b\b\0\b\b\b\b\b\b\b\b\b\b\b\b\b\b\b\0 \\
\0\b\0\0\0\0\0\0\0\0\0\0\0\0\0\b\0\b\0\0\0\0\0\0\0\0\0\0\0\0\0\b\0 \\
\0\b\0\b\b\b\b\b\b\b\b\b\b\b\0\b\0\b\0\b\b\b\b\b\b\b\b\b\b\b\0\b\0 \\
\b\0\b\0\b\0\0\0\0\0\0\0\b\0\b\0\b\0\b\0\b\0\0\0\0\0\0\0\b\0\b\0\b \\
\b\0\b\0\0\0\b\b\b\b\b\0\0\0\b\0\b\0\b\0\0\0\b\b\b\b\b\0\0\0\b\0\b \\
\b\0\b\b\b\b\b\0\0\0\b\b\b\b\b\0\b\0\b\b\b\b\b\0\0\0\b\b\b\b\b\0\b \\
\b\0\0\0\0\0\0\0\b\0\0\0\0\0\0\0\b\0\0\0\0\0\0\0\b\0\0\0\0\0\0\0\b \\
\0\b\b\b\b\b\b\b\b\b\b\b\b\b\b\b\0\b\b\b\b\b\b\b\b\b\b\b\b\b\b\b\0 \\
\0\b\0\0\0\0\0\0\0\0\0\0\0\0\0\b\0\b\0\0\0\0\0\0\0\0\0\0\0\0\0\b\0 \\
\0\b\0\b\b\b\b\b\b\b\b\b\b\b\0\b\0\b\0\b\b\b\b\b\b\b\b\b\b\b\0\b\0 \\
\b\0\b\0\b\0\0\0\0\0\0\0\b\0\b\0\b\0\b\0\b\0\0\0\0\0\0\0\b\0\b\0\b \\
\b\0\b\0\0\0\b\b\b\b\b\0\0\0\b\0\b\0\b\0\0\0\b\b\b\b\b\0\0\0\b\0\b \\
\b\0\b\b\b\b\b\0\0\0\b\b\b\b\b\0\b\0\b\b\b\b\b\0\0\0\b\b\b\b\b\0\b \\
\b\0\0\0\0\0\0\0\b\0\0\0\0\0\0\0\b\0\0\0\0\0\0\0\b\0\0\0\0\0\0\0\b \\
\0\b\b\b\b\b\b\b\b\b\b\b\b\b\b\b\0\b\b\b\b\b\b\b\b\b\b\b\b\b\b\b\0 \\
\0\b\0\0\0\0\0\0\0\0\0\0\0\0\0\b\0\b\0\0\0\0\0\0\0\0\0\0\0\0\0\b\0 \\
\0\b\0\b\b\b\b\b\b\b\b\b\b\b\0\b\0\b\0\b\b\b\b\b\b\b\b\b\b\b\0\b\0 \\
\b\0\b\0\b\0\0\0\0\0\0\0\b\0\b\0\b\0\b\0\b\0\0\0\0\0\0\0\b\0\b\0\b \\
\b\0\b\0\0\0\b\b\b\b\b\0\0\0\b\0\b\0\b\0\0\0\b\b\b\b\b\0\0\0\b\0\b \\
\b\0\b\b\b\b\b\0\0\0\b\b\b\b\b\0\b\0\b\b\b\b\b\0\0\0\b\b\b\b\b\0\b \\
\\ \mbox{ \large $AB^8F'^2F''$ dozens}
}
$
}

\vspace*{5ex}

{\bf Appendix B: Best lower bounds known for $\delta(n)$.}

We exhibit the subsets $S\subseteq\Z^2$ described in Lemma 1 attaining
the largest density known subject to $d(S)\leqs n$.  The resulting
lower bounds on $\delta(n)$ are now known to be sharp for all $n$
except possibly $n=4$, using the Corollary to Prop.~1 for $n\geqs6$,
Prop.~2 for $n\leqs2$, the main Theorem of this paper for $n=3$,
and M.~Moore's work for $n=5$.  See Fig.~5 for further examples of
$S$\/ of density $3/5$ with $d(S)=4$.

\vspace*{5ex}

\centerline{
$
\sq{
\b\0\b\0\b\0\b\0\b\0\b\\
\0\0\0\0\0\0\0\0\0\0\0\\
\b\0\b\0\b\0\b\0\b\0\b\\
\0\0\0\0\0\0\0\0\0\0\0\\
\b\0\b\0\b\0\b\0\b\0\b\\
\0\0\0\0\0\0\0\0\0\0\0\\
\b\0\b\0\b\0\b\0\b\0\b\\
\0\0\0\0\0\0\0\0\0\0\0\\
\b\0\b\0\b\0\b\0\b\0\b\\
\0\0\0\0\0\0\0\0\0\0\0\\
\b\0\b\0\b\0\b\0\b\0\b\\
\\ \delta(0) = 1/4
}
\qquad
\sq{
\b\b\0\b\b\0\b\b\0\b\b\\
\0\0\0\0\0\0\0\0\0\0\0\\
\b\b\0\b\b\0\b\b\0\b\b\\
\0\0\0\0\0\0\0\0\0\0\0\\
\b\b\0\b\b\0\b\b\0\b\b\\
\0\0\0\0\0\0\0\0\0\0\0\\
\b\b\0\b\b\0\b\b\0\b\b\\
\0\0\0\0\0\0\0\0\0\0\0\\
\b\b\0\b\b\0\b\b\0\b\b\\
\0\0\0\0\0\0\0\0\0\0\0\\
\b\b\0\b\b\0\b\b\0\b\b\\
\\ \delta(1) = 1/3
}
\qquad
\sq{
\b\b\b\b\b\b\b\b\b\b\b\\
\0\0\0\0\0\0\0\0\0\0\0\\
\b\b\b\b\b\b\b\b\b\b\b\\
\0\0\0\0\0\0\0\0\0\0\0\\
\b\b\b\b\b\b\b\b\b\b\b\\
\0\0\0\0\0\0\0\0\0\0\0\\
\b\b\b\b\b\b\b\b\b\b\b\\
\0\0\0\0\0\0\0\0\0\0\0\\
\b\b\b\b\b\b\b\b\b\b\b\\
\0\0\0\0\0\0\0\0\0\0\0\\
\b\b\b\b\b\b\b\b\b\b\b\\
\\ \delta(2) = \delta(3) = 1/2
}
$
}

\pagebreak
\vspace*{-10ex}

\centerline{
$
\sq{
\0\0\b\b\b\0\0\b\b\b\0\\
\b\b\0\0\b\b\b\0\0\b\b\\
\0\b\b\b\0\0\b\b\b\0\0\\
\b\0\0\b\b\b\0\0\b\b\b\\
\b\b\b\0\0\b\b\b\0\0\b\\
\0\0\b\b\b\0\0\b\b\b\0\\
\b\b\0\0\b\b\b\0\0\b\b\\
\0\b\b\b\0\0\b\b\b\0\0\\
\b\0\0\b\b\b\0\0\b\b\b\\
\b\b\b\0\0\b\b\b\0\0\b\\
\0\0\b\b\b\0\0\b\b\b\0\\
\\ \delta(4) \geqs 3/5
}
\qquad
\sq{
\0\0\0\0\0\0\0\0\0\0\0\\
\b\b\b\b\b\b\b\b\b\b\b\\
\0\b\0\b\0\b\0\b\0\b\0\\
\0\b\0\b\0\b\0\b\0\b\0\\
\b\b\b\b\b\b\b\b\b\b\b\\
\0\0\0\0\0\0\0\0\0\0\0\\
\b\b\b\b\b\b\b\b\b\b\b\\
\0\b\0\b\0\b\0\b\0\b\0\\
\0\b\0\b\0\b\0\b\0\b\0\\
\b\b\b\b\b\b\b\b\b\b\b\\
\0\0\0\0\0\0\0\0\0\0\0\\
\\ \mbox{alternative $\frac35$ pattern}
}
\qquad
\sq{
\b\b\b\b\b\b\b\0\b\0\0\\
\b\b\0\b\0\0\b\0\b\b\b\\
\0\b\0\b\b\b\b\b\b\b\0\\
\b\b\b\b\b\0\b\0\0\b\0\\
\0\b\0\0\b\0\b\b\b\b\b\\
\0\b\b\b\b\b\b\b\0\b\0\\
\b\b\b\0\b\0\0\b\0\b\b\\
\0\0\b\0\b\b\b\b\b\b\b\\
\b\b\b\b\b\b\0\b\0\0\b\\
\b\0\b\0\0\b\0\b\b\b\b\\
\b\0\b\b\b\b\b\b\b\0\b\\
\\ \delta(5) = 9/13\ \cite{Moore}
}
$
}

\vspace*{3ex}

\centerline{
$
\sq{
\b\b\0\b\b\b\b\0\b\b\\
\b\b\b\b\0\b\b\b\b\0\\
\b\0\b\b\b\b\0\b\b\b\\
\b\b\b\0\b\b\b\b\0\b\\
\0\b\b\b\b\0\b\b\b\b\\
\b\b\0\b\b\b\b\0\b\b\\
\b\b\b\b\0\b\b\b\b\0\\
\b\0\b\b\b\b\0\b\b\b\\
\b\b\b\0\b\b\b\b\0\b\\
\\ \delta(6) = 4/5
}
\qquad
\sq{
\b\b\b\b\b\b\b\b\b\\
\b\0\b\b\0\b\b\0\b\\
\b\b\b\b\b\b\b\b\b\\
\b\b\b\b\b\b\b\b\b\\
\b\0\b\b\0\b\b\0\b\\
\b\b\b\b\b\b\b\b\b\\
\b\b\b\b\b\b\b\b\b\\
\b\0\b\b\0\b\b\0\b\\
\b\b\b\b\b\b\b\b\b\\
\\ \delta(7) = 8/9
}
\qquad
\sq{
\b\b\b\b\b\b\b\b\b\\
\b\b\b\b\b\b\b\b\b\\
\b\b\b\b\b\b\b\b\b\\
\b\b\b\b\b\b\b\b\b\\
\b\b\b\b\b\b\b\b\b\\
\b\b\b\b\b\b\b\b\b\\
\b\b\b\b\b\b\b\b\b\\
\b\b\b\b\b\b\b\b\b\\
\b\b\b\b\b\b\b\b\b\\
\\ \delta(8) = 1
}
$
}

\vspace*{5ex}

{\bf Appendix C: Patterns attaining $\delta_4(n)$ and $\delta_6(n)$.}

We exhibit the patterns described in Prop.~3 attaining maximal
densities for the 4- and 6-point neighborhoods.

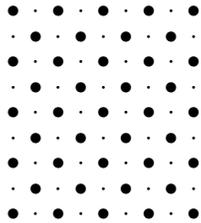
\vspace*{3ex}

\centerline{
$
\sq{
\b\0\b\0\b\0\b\0\b\\
\0\b\0\b\0\b\0\b\0\\
\b\0\b\0\b\0\b\0\b\\
\0\b\0\b\0\b\0\b\0\\
\b\0\b\0\b\0\b\0\b\\
\0\b\0\b\0\b\0\b\0\\
\b\0\b\0\b\0\b\0\b\\
\0\b\0\b\0\b\0\b\0\\
\b\0\b\0\b\0\b\0\b\\
\\ \delta_4(0) = \delta_4(1) = 1/2
}
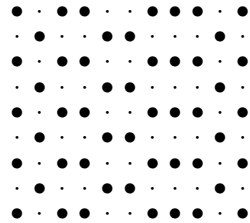
\qquad\quad
\sq{
\b\0\b\b\0\0\b\b\b\0\b\\
\0\b\0\0\b\b\0\0\0\b\0\\
\b\0\b\b\0\0\b\b\b\0\b\\
\0\b\0\0\b\b\0\0\0\b\0\\
\b\0\b\b\0\0\b\b\b\0\b\\
\0\b\0\0\b\b\0\0\0\b\0\\
\b\0\b\b\0\0\b\b\b\0\b\\
\0\b\0\0\b\b\0\0\0\b\0\\
\b\0\b\b\0\0\b\b\b\0\b\\
\\ \delta_4(0) = 1/2\mbox{\ ``chicken\ wire''}
}
$
}

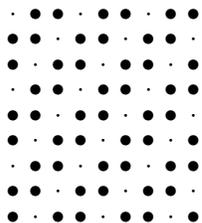
\vspace*{3ex}

\centerline{
$
\sq{
\0\b\b\0\b\b\0\b\b\\
\b\b\0\b\b\0\b\b\0\\
\b\0\b\b\0\b\b\0\b\\
\0\b\b\0\b\b\0\b\b\\
\b\b\0\b\b\0\b\b\0\\
\b\0\b\b\0\b\b\0\b\\
\0\b\b\0\b\b\0\b\b\\
\b\b\0\b\b\0\b\b\0\\
\b\0\b\b\0\b\b\0\b\\
\\ \delta_4(2) = 2/3
}
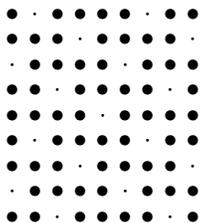
\qquad
\sq{
\b\0\b\b\b\b\0\b\b\\
\b\b\b\0\b\b\b\b\0\\
\0\b\b\b\b\0\b\b\b\\
\b\b\0\b\b\b\b\0\b\\
\b\b\b\b\0\b\b\b\b\\
\b\0\b\b\b\b\0\b\b\\
\b\b\b\0\b\b\b\b\0\\
\0\b\b\b\b\0\b\b\b\\
\b\b\0\b\b\b\b\0\b\\
\\ \delta_4(3) = 4/5
}
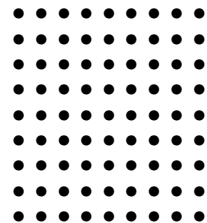
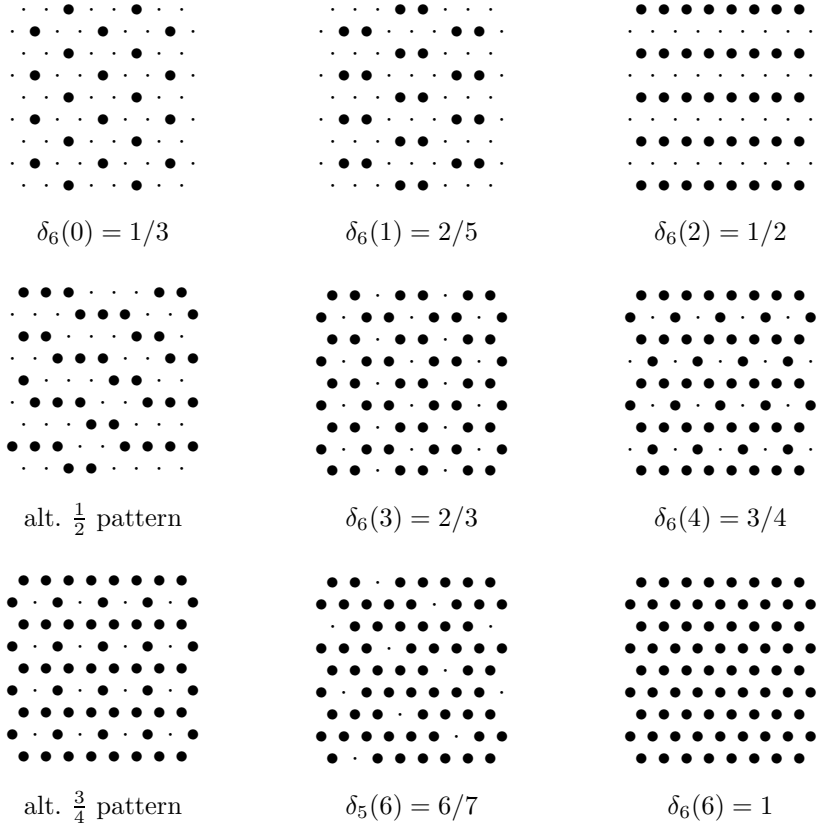
\qquad
\sq{
\b\b\b\b\b\b\b\b\b\\
\b\b\b\b\b\b\b\b\b\\
\b\b\b\b\b\b\b\b\b\\
\b\b\b\b\b\b\b\b\b\\
\b\b\b\b\b\b\b\b\b\\
\b\b\b\b\b\b\b\b\b\\
\b\b\b\b\b\b\b\b\b\\
\b\b\b\b\b\b\b\b\b\\
\b\b\b\b\b\b\b\b\b\\
\\ \delta_4(4) = 1
}
$
}

\pagebreak

\renewcommand{\arraystretch}{.693}
\centerline{
$
\sq{
 \0\0\b\0\0\b\0\0 \\
\0\b\0\0\b\0\0\b\0\\
 \0\0\b\0\0\b\0\0 \\
\0\b\0\0\b\0\0\b\0\\
 \0\0\b\0\0\b\0\0 \\
\0\b\0\0\b\0\0\b\0\\
 \0\0\b\0\0\b\0\0 \\
\0\b\0\0\b\0\0\b\0\\
 \0\0\b\0\0\b\0\0 \\
\\ \delta_6(0)=1/3
}
\quad\qquad
\sq{
 \0\0\0\b\b\0\0\0 \\
\0\b\b\0\0\0\b\b\0\\
 \0\0\0\b\b\0\0\0 \\
\0\b\b\0\0\0\b\b\0\\
 \0\0\0\b\b\0\0\0 \\
\0\b\b\0\0\0\b\b\0\\
 \0\0\0\b\b\0\0\0 \\
\0\b\b\0\0\0\b\b\0\\
 \0\0\0\b\b\0\0\0 \\
\\ \delta_6(1)=2/5
}
\quad\qquad
\sq{
 \b\b\b\b\b\b\b\b \\
\0\0\0\0\0\0\0\0\0\\
 \b\b\b\b\b\b\b\b \\
\0\0\0\0\0\0\0\0\0\\
 \b\b\b\b\b\b\b\b \\
\0\0\0\0\0\0\0\0\0\\
 \b\b\b\b\b\b\b\b \\
\0\0\0\0\0\0\0\0\0\\
 \b\b\b\b\b\b\b\b \\
\\ \delta_6(2)=1/2
}
$
}

\vspace*{3ex}

\centerline{
$
\sq{
 \b\b\b\0\0\0\b\b \\
\0\0\0\b\b\b\0\0\b\\
 \b\b\0\0\0\b\b\0 \\
\0\0\b\b\b\0\0\b\b\\
 \b\0\0\0\b\b\0\0 \\
\0\b\b\b\0\0\b\b\b\\
 \0\0\0\b\b\0\0\0 \\
\b\b\b\0\0\b\b\b\b\\
 \0\0\b\b\0\0\0\0 \\
\\ \mbox{alt.\ $\frac12$ pattern}
}
\quad\qquad
\sq{
 \b\b\0\b\b\0\b\b \\
\b\0\b\b\0\b\b\0\b\\
 \b\b\0\b\b\0\b\b \\
\b\0\b\b\0\b\b\0\b\\
 \b\b\0\b\b\0\b\b \\
\b\0\b\b\0\b\b\0\b\\
 \b\b\0\b\b\0\b\b \\
\b\0\b\b\0\b\b\0\b\\
 \b\b\0\b\b\0\b\b \\
\\ \delta_6(3)=2/3
}
\quad\qquad
\sq{
 \b\b\b\b\b\b\b\b \\
\b\0\b\0\b\0\b\0\b\\
 \b\b\b\b\b\b\b\b \\
\0\b\0\b\0\b\0\b\0\\
 \b\b\b\b\b\b\b\b \\
\b\0\b\0\b\0\b\0\b\\
 \b\b\b\b\b\b\b\b \\
\0\b\0\b\0\b\0\b\0\\
 \b\b\b\b\b\b\b\b \\
\\ \delta_6(4)=3/4
}
$
}

\vspace*{3ex}

\centerline{
$
\sq{
 \b\b\b\b\b\b\b\b \\
\b\0\b\0\b\0\b\0\b\\
 \b\b\b\b\b\b\b\b \\
\b\0\b\0\b\0\b\0\b\\
 \b\b\b\b\b\b\b\b \\
\b\0\b\0\b\0\b\0\b\\
 \b\b\b\b\b\b\b\b \\
\b\0\b\0\b\0\b\0\b\\
 \b\b\b\b\b\b\b\b \\
\\ \mbox{alt.\ $\frac34$ pattern}
}
\quad\qquad
\sq{
 \b\b\0\b\b\b\b\b \\
\b\b\b\b\b\0\b\b\b\\
 \0\b\b\b\b\b\b\0 \\
\b\b\b\0\b\b\b\b\b\\
 \b\b\b\b\b\0\b\b \\
\b\0\b\b\b\b\b\b\0\\
 \b\b\b\0\b\b\b\b \\
\b\b\b\b\b\b\0\b\b\\
 \b\0\b\b\b\b\b\b \\
\\ \delta_5(6)=6/7
}
\quad\qquad
\sq{
 \b\b\b\b\b\b\b\b \\
\b\b\b\b\b\b\b\b\b\\
 \b\b\b\b\b\b\b\b \\
\b\b\b\b\b\b\b\b\b\\
 \b\b\b\b\b\b\b\b \\
\b\b\b\b\b\b\b\b\b\\
 \b\b\b\b\b\b\b\b \\
\b\b\b\b\b\b\b\b\b\\
 \b\b\b\b\b\b\b\b \\
\\ \delta_6(6)=1
}
$
}


\begin{thebibliography}{0}
\bibitem{WW} Berlekamp, E.R., Conway, J.H., Guy, R.K.:
 {\em Winning Ways For Your Mathematical Plays, II:
 Games In Particular}.  London: Academic Press, 1982.
\bibitem{BC} Buckingham, D., Callahan, P.: Tight bounds
 on periodic cell configurations in Life.  Manuscript, Baltimore 1997.
\bibitem{Greg} Kuperberg, G.:  e-mail communication 1994, 1997.
\bibitem{Moore} Moore, M.: Posting on {\tt comp.theory.cell-automata},
 Dec.~95 (archived at {\sf http://www.krl.caltech.edu/%
 $\sim$brown/news/cell-automata-html/0644.html});
 \hbox{e-mail} communication 1996.
\bibitem{Robinson} Robinson, R.M.: Undecidability and nonperiodicity
 for tilings of the plane.  {\em Invent.\ Math.}\ {\bf 12} (1971),
 177--209.
\bibitem{SK} Saaty, T.L., Kainen, P.C.: {\em The Four-Color Problem:
 Assaults amd Conquest}.  New York: McGraw-Hill, 1977.
\bibitem{TECC} MacWilliams, F.J., Sloane, N.J.A.:
{\em The Theory of Error-Correcting Codes.}
Amsterdam: North-Holland, 1977.
\end{thebibliography}
\end{document}